# The geometry of Brownian surfaces[*]

**Rémi Léandre**

*Institut de Mathématiques, Université de Bourgogne, 21000, Dijon, France*
*e-mail:* `remi.leandre@u-bourgogne.fr`

**Abstract:** Motivated by Segal's axiom of conformal field theory, we do a survey on geometrical random fields. We do a history of continuous random fields in order to arrive at a field theoretical analog of Klauder's quantization in Hamiltonian quantum mechanic by using infinite dimensional Airault-Malliavin Brownian motion.

## 1. Introduction

Let us consider $m$ smooth vector fields on $R^d$, $X_i$. Suppose that we are in an uniform elliptic situation, that is the following quadratic form

$$\xi \to \sum <X_i(x), \xi>^2 \tag{1}$$

is uniformly invertible on $R^d$. In such a case, we can introduce the inverse quadratic form $g(x)$ and a measure $dg = \det(g(x))^{-1/2} dx$. These data are invariant under change of coordinates on $R^d$. A 1-form can be assimilated via the Riemannian metric $g(x)$ to a vector field. We have the following integration by parts formula for any vector field on $R^d$ with compact support:

$$\int_{R^d} Xf \, dg = \int_{R^d} f \operatorname{div} X \, dg \tag{2}$$

Let $\operatorname{grad} f$ the vector field associated to $f$ via the Riemannian metric. The Laplace Beltrami operator is

$$\Delta = \operatorname{div} \operatorname{grad} \tag{3}$$

Since these data are compatible with change of coordinates, we can consider via local charts a manifold $M$ endowed with a Riemannian structure. The notion of Riemannian measure is intrinsic and the Laplace-Beltrami operator $\Delta$ is intrinsically associated to the Riemannian manifold $M$.

Let us consider m independent Brownian motions $B_i$ on $R$ and the Stratonovitch differential equation:

$$dx_t(x) = \sum X_i(x_t(x))dB_t^i; \quad x_0(x) = x \tag{4}$$

---

[*]This is an original survey paper

It generates a Markov semi-group whose generator is $1/2 \sum X_i^2$. Moreover, there exists a vector field such that:

$$-\sum X_i^2 = \Delta + X_0 \tag{5}$$

Since the Itô formula in Stratonovitch sense is the classical one, (4) is independent of the system of coordinates chosen, and we can consider (4) on a manifold.

For that, let us recall quickly the theory of stochastic processes on a manifold. Let $M$ be a compact manifold. Let $t \to x_t$ be a random continuous process on the manifold adapted to a filtration $F_t$ endowed with a Probability measure $P$. $t \to x_t$ is said to be a semimartingale with values in $M$ if for all smooth function $f$ on $M$, $t \to f(x_t)$ is a real semimartingale. Let us consider some smooth vector fields $X_i$ on $M$. Alternatively, they can be considered as first order differential operators on the space of smooth functions on $M$ satisfying the Leibniz rule. Or they can be considered as smooth of the tangent bundle on $M$. The solution of the Stratonovitch differential equation

$$dx_t(x) = \sum_{i>0} X_i(x_t(x))dB_t^i + X_0(x_t(x))dt; \quad x_0(x) = x \tag{6}$$

is characterized by the following condition:

$$df(x_t(x)) = \sum_{i>0} X_i f(x_t(x))dB_t^i + X_0 f(x_t(x))dt \tag{7}$$

for each smooth function $f$, where the vector fields $X_i$ are considered as differential operators. Moreover, if the manifold is compact, the solution of the Stratonovitch differential equation has a unique solution which is a continuous semi-martingale.

If $t \to x_t$ is a continuous semi-martingale and $\omega$ is a smooth 1-form, that is a smooth section of the cotangent bundle of $M$, we can define the stochastic Stratonovitch integral $\int_0^1 <\omega(x_t), dx_t>$ as limit in $L^2$ of the classical random integral

$$\int_0^1 <\omega(x_t^n), dx_t^n> \tag{8}$$

where $t \to x_t^n$ is a suitable polygonal approximation of $t \to x_t$. If $f$ is a smooth real function on $M$, $df$ is a 1-form and the Itô-Stratonovitch formula says that almost surely:

$$f(x_t) - f(x_s) = \int_s^t <df(x_u), dx_u> \tag{9}$$

for $s < t$.

Solution of Stratonovitch differential equations realize measure on the *continuous* path space $P(M)$ constituted of continuous paths from $[0,1]$ into $M$. This realizes random paths on a manifold. There are two geometries involved:

– The source geometry is the segment $[0, 1]$ (or the circle, if we conditionate by $x_1(x) = x$, that is if we consider the Bridge of the diffusion $t \to x_t(x)$ or a tree if we consider branching processes).

– The target geometry constituted of the Riemannian manifold $M$.

If we imbed isometrically the compact manifold $M$ into $R^m$ (It is possible by Nash theorem), we can consider the orthogonal projection $\Pi(x)$ from $R^m$ into $T_x(M)$, the tangent space of $M$ at x. We extend $\Pi$ for $x \in R^m$. We can consider the diffusion:

$$dx_t(x) = \Pi(x_t(x))dB_t; \quad x_0(x) = x \tag{10}$$

where $x$ belongs to $M$ and $B_t$ is a Brownian motion on $R^m$. Since the Itô formula in Stratonovitch Calculus is the traditional one, $t \to x_t(x)$ is a process on $M$ which realizes the Brownian motion on $M$. It is the Markov process whose generator is the Laplace-Beltrami operator on $M$. The source geometry is very simple:

– It is $[0,1]$ if we consider the diffusion. $[0,1]$ is endowed with its canonical Riemannian structure, but we could choose another Riemannian structure by looking at the equation $dx_t(x) = h(t)\Pi(x_t(x))dB_t$.
– It is the circle if we consider the Brownian bridge.

In theoretical physics, people look at random fields where the source is a Riemannian manifold $\Sigma$ with possible boundaries and whose target spaces are manifolds (see [46, 47, 48] for surveys). There are two geometries involved in two dimensional field theories:

– The geometry of the world sheet, a Riemannian surface with boundary.
– The geometry of the target manifold.

If we consider the case where the world sheet is a square $[0,1] \times [0,1]$ or a cylinder $[0,1] \times S^1$, the random field $(t,s) \to x_{t,s}$ realizes an infinite dimensional process $t \to (s \to x_{t,s})$ with values in the path space of the manifold or in the loop space of the manifold. We say that we are in presence of a $1+1$ dimensional field theory, the first 1 describing the time of the dynamic and the second 1 denoting the dimension of the internal time of the state space. If we consider a more complicated Riemannian surface, the loops are interacting: we refer to the survey of Mandelstam [102] or Witten [126] about that.

There are two possible theories:

– We consider open strings, that is possibly interacting processes on the path space on a manifold.
– We consider closed strings, that is possibly interacting processes on the loop space.

If we consider processes on the loop space, the axioms of Segal [119] of conformal field theory are the followings:

**Segal's axiom:** Consider the set of possibly disconnected Riemannian surfaces $\Sigma$ with oriented parametrization $p_i, i \in I$ of the boundary loops, negative (positive) for $i \in I_-(I_+)$ and with a Riemannian structure $g$, trivial around the boundary. Let us consider an Hilbert space $H$ with an antiunitary involution $P$.

A real conformal field theory is an assignment

$$(\Sigma, p_i, g) \to A(\Sigma, p_i, g) \tag{11}$$

where
$$A(\Sigma, p_i, g) : \otimes_{i \in I_-} H \to \otimes_{i \in I_+} H \tag{12}$$

are trace-classes operators (empty tensor are equal to $C$) satisfying to the following properties:

*Property 1:* If $(\Sigma, p_i, g)$ is the disjoint union of $(\Sigma^\alpha, p_{i_\alpha}^\alpha, g^\alpha)$, then:
$$A(\Sigma, p_i, g) = \otimes_\alpha A(\Sigma^\alpha, p_{i_\alpha}^\alpha, g^\alpha) \tag{13}$$

*Property 2:* If we reverse the sense of the time of $p_{i_0}$, we get another loop called $\tilde{p}_{i_0}$ such that if $i_0 \in I_-$:
$$< A(\Sigma, \tilde{p}_{i_0}, p_i, g) x_{i_0} \otimes x, y > = < A(\Sigma, p_i, g) x, P x_{i_0} \otimes y > \tag{14}$$

*Property 3:* If $F$ is a conformal diffeomorphism from $\Sigma^1$ into $\Sigma^2$, then:
$$A(\Sigma^1, p_i^1, F^* g^2) = A(\Sigma^2, F \circ p_i^1, g^2) \tag{15}$$

*Property 4:* If $\Sigma'$ is constructed from $\Sigma$ by identifying the boundary loops $i_1 \in I_-$ and $i_2 \in I_+$, we have:
$$A(\Sigma', p_{i'}, g) = \text{Tr}_{i_1, i_2} A(\Sigma, p_i, g) \tag{16}$$

where $i'$ is different from $i_1$ and $i_2$ and $\text{Tr}_{i_1, i_2}$ is the trace between factors $i_1$ and $i_2$ in the tensor products of $H$.

*Property 5:* If $\tilde{\Sigma}$ is the complex conjugate of $\Sigma$, then:
$$A(\tilde{\Sigma}, p_i, g) = A(\Sigma, p_i, g)^* \tag{17}$$

*Property 6:* If $\sigma$ is a real smooth function on $\Sigma$ vanishing in a boundary of the Riemann surface $\Sigma$, then:
$$A(\Sigma, p_i, \exp[\sigma] g) = \exp[\frac{ci}{2\pi} \int_\Sigma (1/2 \partial \sigma \wedge \overline{\partial} \sigma + R_g \sigma)] A(\Sigma, p_i, g) \tag{18}$$

where $R_g$ is the curvature form of the metric $g$ on the surface $\Sigma$ and $c$ is a constant.

*Property 6* is called the conformal anomaly. Namely, these axioms as it was noticed by Gawedzki (Ref. [46], pp. 106–107) may be deduced from the physicist intuitive representative of the amplitude $A$ by formal functional integrals. $H$ is a space of function over the loop space $L(M)$ of the target space $M$. $P$ arises on the time reversal on the loop $L(M)$ combined with the complex conjugation. The amplitudes $A$ are represented as formal integrals on maps $x_. : \Sigma \to M$ fixed on the boundary of $\Sigma$ (This means we consider a kind of infinite dimensional Brownian bridge):
$$\int_{x_. \circ p_i = f_i} d\mu(x_.) = \int_{x_. \circ p_i = f_i} \exp[-I(x_.)] dD(x_.) = A(\Sigma, p_i, g)(f_i) \tag{19}$$

where $dD(x_.)$ is the formal Lebesgue measure over the set of maps $x_.$ and $I(x_.)$ the energy of $x_.$:

$$I(x_.) = \int_\Sigma <Dx_S, Dx_S>_{x_S} dg(S) \qquad (20)$$

where $S$ denotes the generic element of the Riemannian surface $\Sigma$. The quantities (19) lead to infinities, which are regularized by ad-hoc procedures in physics. Let us consider the case where $M = R$ endowed with the constant metric, and in order to simplify, the case where we have no boundary in the Riemannian surface $\Sigma$. Let us consider the partition function of the theory:

$$Z = \int_{x_.} d\mu(x_.) \qquad (21)$$

and the Gaussian measure $Z^{-1}d\mu(x_.)$. The random field $x_.$ is Gaussian and is called the free field. Let us recall some basic backgrounds on Gaussian random fields parametrized by $\Sigma$. We suppose in order to simplify that the Riemann surface $\Sigma$ has no boundary. Let $\Delta_\Sigma$ be the Laplacian on $\Sigma$. Let be the action

$$I(x(.)) = \int_\Sigma <(\Delta_\Sigma + I)x(.), x(.)> dg \qquad (22)$$

Let $x_n(s)$ be the normalized eigenfunctions associated to the eigenvalues $\lambda_n$ of $\Delta_\Sigma + I$. Let us remark that $\lambda_n > 0$ for all $n$ due to $I$. The Gaussian field associated to $I(x(.))$ can be represented as

$$x(.) = \sum \frac{1}{\lambda_n} B_n x_n(.) \qquad (23)$$

where $B_n$ are independent centered normalized Gaussian variables. There are some problems to know the space where the previous series converges. If we consider a smooth function $f$ on $M$, it can be represented by Sobolev imbedding theorem as

$$f(.) = \sum_n \alpha_n x_n(.) \qquad (24)$$

where the sequence $\alpha_n$ is quickly decreasing. This shows us that

$$\int f(S)x(S)dg(S) = \sum \frac{\alpha_n}{\lambda_n} B_n \qquad (25)$$

is a convergent series of random variables. This shows us that the random field can be defined as a random distribution.

In order to show that the formal random field defined by (22) is not a true random field, we can compute formally by using (23) $E[x(S)x(S')]$. It is given by

$$\sum \frac{1}{\lambda_n} x_n(S) x_n(S') \qquad (26)$$

We recognize in the previous expression the kernel of $(\Delta_\Sigma + I)^{-1}$ (or the Green kernel associated to $\Delta_\Sigma + I$), which has a logarithmic singularity when $S \to S'$ $\log d(S, S')$ where $d$ is the Riemannian distance on $\Sigma$.

In our situation, we have no mass term $I$ and the treatment is a little big more complicated due to the presence of zero modes. Let us consider the Laplace-Beltrami operator on $M$ and the associated Green kernel $G(S, S')$. We have $E[x_S x_{S'}] = G(S, S')$. We still have $G(S, S) = \infty$. This explain, following Nelson [108], that the random field $x_{\cdot}$ is a generalized process: it is a random distribution [108] and we consider smeared random fields. It is difficult to state what are random distributions with values in a manifold.

But there are simpler random fields than the free field: for instance, the Brownian sheet is a true random field. The goal of this survey is to recall them and to describe the history in order to arrive to a theory of (continuous!) random fields parametrized by any Riemannian surface with values in any target space.

For surveys about the physicists models, we refer to the survey of Witten [126] and Gawedzki [46, 47, 48]. For the relation between analysis on loop space and mathematical physics, we refer to the surveys of Albeverio [3] and Léandre [80, 83].

## 2. Cairoli equation and Brownian sheet

In this part, we are concerned by a world-sheet which is a square $[0, 1] \times [0, 1]$. Namely, in this case, there is an order on the world-sheet such that we can apply martingale theory in order to study stochastic differential equations.

Let $S = (t, s) \in [0, 1]^2$. We consider the white noise $\eta(S)$. It is a generalized Gaussian process with average 0 and formal covariance given by

$$E[\eta(S)\eta(S')] = \delta_S(S') \qquad (27)$$

where $\delta_S(S')$ is the Dirac mass in $S$. We can consider a measurable set $O$ of the square and we can define

$$B(O) = \int_O \eta(S) dS \qquad (28)$$

$B(O)$ is a Gaussian variable of average 0. Moreover,

$$E[B(O)B(O')] = m(O \cap O') \qquad (29)$$

where $m$ is the Lebesgue measure on the square.

On the square, there is a natural order. If $S = (t, s)$, we denote by $[0, S] = [0, t] \times [0, s]$ and we can introduce the Brownian sheet:

$$B(S) = \int_{[0,S]} \eta(S') dS' \qquad (30)$$

It is a continuous process. Namely, we can use (29) and the fact that $B(S) - B(S')$ is a Gaussian random variable in order to find a positive $\alpha$ such that:

$$E[|B(S) - B(S')|^p]^{1/p} \leq Cd(S, S')^\alpha \qquad (31)$$

The result arises by Kolmogorov lemma.

Let $F(S)$ be the $\sigma$-algebra spanned by $B(S')$ for $S' \leq S$. If $S_1 > S$, we get almost surely

$$B(S) = E[B(S_1)|F(S)] \tag{32}$$

We say that $s \to B(S)$ is a martingale with respect to the filtration $F(S)$.

Let $S \to h(S)$ be a continuous bounded process such that $h(S)$ is $F(S)$-measurable. We can define the Itô stochastic integral

$$V(S) = \int_{[0,S]} h(S')\delta B(S') \tag{33}$$

$S \to V(S)$ is still a martingale and Itô's isometry formula states that:

$$E[V(S)^2] = E[\int_{[0,S]} h^2(S')dS'] \tag{34}$$

We refer to [30] for a theory of Itô integral for multiparameter random processes.

This equality, by using a generalization in the two parameter context of the Peano approximation of a diffusion, allows to consider Cairoli equation [22]: let be $m+1$ vector fields $X_i$ on $R^d$ with bounded derivatives of all orders and $m$ independent Brownian sheets $B^i(S)$. We consider the two-parameter Itô equation:

$$\delta x(S) = \sum X_i(x(S))\delta B^i(S) + X_0(S)dS; \quad x(0) = x \tag{35}$$

This means that:

$$x(S) = x + \int_{[0,S]} X_i(x(S'))\delta B^i(S') + \int_{[0,S]} X_0(x(S'))dS' \tag{36}$$

**Theorem 1** *(Cairoli) Equation (35) has a unique solution.*

Let us introduce the Hilbert space $H_1$ of maps from the square into $R^m$ $h$ such that $h(S) = \int_{[0,S]} k(S')dS'$ endowed with the Hilbert structure $\|h\|^2 = \int_{[0,1]^2} |k(S')|^2 dS'$. We remark that:

$$\|h\|^2 = \int_{[0,1]^2} |\frac{\partial^2}{\partial s \partial t} h(S)|^2 ds \tag{37}$$

(compare with (22)).

Formally, the Brownian sheet is the Gaussian measure on $H_1$ defined, in a heuristic physicist way by:

$$Z^{-1}\exp[-\|h\|^2/2]dD(h) \tag{38}$$

where $dD(h)$ is the heuristic Lebesgue measure on $H_1$. In fact, unlike the free field, which is a generalized process, this measure lives on the space of continuous functions from the square with values in $R^m$ endowed with the uniform norm $\|.\|_\infty$. The link between the physicist heuristic point of view and the rigorous probabilistic point of view is expressed by large deviations results [11]. Let us consider $\epsilon B(.)$ when $\epsilon \to 0$. It is given by the formal measure

$$Z_\epsilon^{-1} \exp[-\frac{\|h\|^2}{2\epsilon^2}]dD(h) \tag{39}$$

Let $A$ be a Borelian subset on the space of continuous functions from the square into $R^m$, let Int$A$ be its interior for the uniform topology and Clos $A$ its closure for the uniform topology. (39) explains heuristically that:

$$-\inf_{h \in \text{Int } A} \|h\|^2 \leq \underline{\lim}_{\epsilon \to 0} 2\epsilon^2 \log P\{\epsilon B(.) \in A\} \tag{40}$$

$$\overline{\lim}_{\epsilon \to 0} 2\epsilon^2 \log P\{\epsilon B(.) \in A\} \leq -\inf_{h \in \text{Clos } A} \|h\|^2 \tag{41}$$

Let us consider the equation

$$dx(S)(\epsilon) = \epsilon \sum X_i(x_S(\epsilon))\delta B^i(S); \quad x(0)(\epsilon) = x \tag{42}$$

and its skeleton:

$$dx(S)(h) = \sum X_i(x(S)(h))k(S)^i dS; \quad x(0)(h) = x \tag{43}$$

If $\|h\|$ is bounded, the set of solution of (43) is relatively compact in the set of continuous functions from the square into $R^d$ endowed with the uniform topology. This is due to Ascoli theorem.

Doss and Dozzi [29] have shown the following quasi-continuity lemma:

**Lemma 2** *For all $a, R, \rho > 0$, there exists $\alpha, \epsilon_0$ such that:*

$$P\{\|x(.)(\epsilon) - x(.)(h)\|_\infty > \rho; \|\epsilon B(.) - h(.)\|_\infty < \alpha\} \leq \exp[-R/\epsilon^2] \tag{44}$$

*for $\epsilon \in ]0, \epsilon_0]$ and $h$ such that $\|h\|^2 \leq a$.*

The proof of this lemma is based upon a generalization in the multiparameter context of the exponential inequality for martingales and of the Cameron-Martin formula.

If $A$ is a Borelian of the set of continuous functions from the square into, $R^d$, we put $\Lambda(A) = \inf_{x(.)(h) \in A} \|h\|^2$.

Doss-Dozzi [29] deduced the following large deviation estimates for the two-parameter diffusion $x(.)(\epsilon)$.

**Theorem 3** *(Wentzel-Freidlin estimates): When $\epsilon \to 0$:*

$$-\Lambda(\text{Int } A) \leq \underline{\lim} 2\epsilon^2 \log P\{x(.)(\epsilon) \in A\} \tag{45}$$

$$\overline{\lim} 2\epsilon^2 P\{x(.)(\epsilon) \in A\} \leq -\Lambda(\text{Clos } A) \tag{46}$$

This shows that in order to estimate $P\{x(.)(\epsilon) \in A\}$, we can replace formally $\epsilon B(.)$ by $h$ with the formal measure (38).

We refer to [16, 114] for improvements of this theorem.

Let us suppose that

$$C'|\xi|^2 \geq \sum_{i>0} < X_i(x), \xi >^2 \geq C|\xi|^2 \tag{47}$$

for some $C' > C > 0$.

Under this ellipticity hypothesis, Nualart and Sanz [113] have shown by using Malliavin Calculus:

**Theorem 4** *The law of $x(S)(\epsilon)$ has a smooth density $q_\epsilon(S)(x,y)$ with respect of the Lebesgue measure on $R^d$ provided that $st \neq 0$.*

Let us put

$$d_S^2(x,y) = \inf_{x(S)(h)=y} \|h\|^2 \tag{48}$$

Under (47), $d_S^2(x,y)$ is finite and continuous in $x$ and $y$. Léandre and Russo [99] have mixed Malliavin Calculus and Wentzel-Freidlin estimates in order to show:

**Theorem 5** *(Varadhan estimates): We have if $st \neq 0$ when $\epsilon \to 0$:*

$$\lim 2\epsilon^2 \log q_\epsilon(S)(x,y) = -d_S^2(x,y) \tag{49}$$

The problem of this theory is that we consider Itô equation. This leads to some problems if we want to constrain the two parameter diffusion to live over a manifold. Norris [110] extending a previous work of Hajek [52] has defined a two-parameter Stratonovitch Calculus which allows him to solve this problem. Namely, Itô formula in two parameters Itô Calculus has no geometrical meaning. Stratonovitch Calculus has a geometrical meaning. This allows Eells-Elworthy-Malliavin to get an intrinsic construction of the Brownian motion on a compact manifold. Norris equation is motivated by a multiparameter extension of the construction of Eells-Elworthy-Malliavin of the Brownian motion on a manifold. This requires the introduction of some geometrical definitions.

Let $M$ be a Riemannian manifold. It inherits the Levi-Civita connection $\nabla$. It is characterized by the fact that it is a metric connection without torsion. If $X$ is a vector field and $Y$ a vector fields, $\nabla_Y X$ is still a vector field. Moreover,

$$Y < X^1, X^2 > = < \nabla_Y X^1, X^2 > + < X^1, \nabla_Y X^2 > \tag{50}$$

(The Levi-Civita connection preserves the metric) and

$$\nabla_Y X - \nabla_X Y = [X, Y] \tag{51}$$

(The Levi-Civita connection is without torsion). $\nabla_Y X$ is tensorial in $Y$ and a first order operator in $X$. This means, if $f$ is a smooth function on $M$

$$\nabla_Y (fX) = Y(f)X + f\nabla_Y X \tag{52}$$

In local coordinates,
$$\nabla_Y X = Y(X) + \Gamma_Y X \tag{53}$$
where $\Gamma$ denotes the Christoffel symbol matrix.

Is $s \to x_s$ is a $C^1$ paths with values in the manifold and $s \to X_s$ a path in the tangent bundle over $x_s$, we can define in a intrinsic way:
$$\nabla_s X_s = d/ds X_s + \Gamma_{d/ds x_s} X_s \tag{54}$$

The Levi-Civita connection pass to the cotangent bundle by duality. Moreover, let $O(M)$ be the $O(d)$ principal bundle on $M$ constituted of isometries from $R^d$ into $T_x(M)$. The Levi-Civita connection pass to $O(M)$.

Christoffel symbols and (54) have only a local meaning. In order to perform the computations, we have to see how these quantities behave under local change of coordinates. There is a way to avoid these technicalities. Namely, we can show that all linear bundle $E$ endowed with a linear connection $\nabla^E$ is a subbundle of a trivial bundle endowed with the projection connection. Therefore, Christoffel symbols become globally defined.

Eells-Elworthy-Malliavin equation on a Riemannian manifold is
$$dx_t = \tau_t dB_t \tag{55}$$
where $B_t$ is a Brownian motion in the tangent space of $M$ at $x$ and where $\tau_t$ is the parallel transport for the Levi-Civita connection. Parallel transport after performing the previous trivialization is given by
$$d\tau_t = -\Gamma_{dx_t} \tau_t \tag{56}$$
$\tau_t$ realizes a isometry from the tangent space at $x$ to the tangent space at $x_t$.

We can generalize this notion: let $s \to y_s$ be a continuous semi-martingale in $M$ and $E$ a linear bundle on $M$ endowed with a connection $\nabla^E$. The parallel transport along the random path $y_t$ is given by the equation:
$$d\tau_t^E = -\Gamma_{dy_t}^E \tau_t^E \tag{57}$$
If we consider a path $X_t = \tau_t^E h_t$ over $y_t$, we have
$$\nabla_t^E \tau_t^E h_t = \tau_t^E d/dt h_t \tag{58}$$
where $t \to H_t$ is a finite energy path in the fiber of $E$ over the starting point of the leading semi-martingale.

After trivializing globally the tangent bundle for the Levi-Civita connection by adding a suitable auxiliary bundle, Eells-Elworthy-Malliavin equation becomes more tractable: it is
$$dx_t = \tau_t dB_t; \quad d\tau_t = -\Gamma_{\tau_t dB_t} \tau_t \tag{59}$$

Let us consider a semi-martingale $V(S)$ defined by (33). $s \to V(t,s)$ is a semi-martingale and $t \to V(t,s)$ is a semi-martingale. We can define the Stratonovitch differential $d_s V(S)$ and $d_t V(S)$ of these semi-martingales. Norris ([110], p. 292) defines too the Stratonovitch differential $d_t d_s V(S)$. These Stratonovitch differentials satisfy to:

$$d_s f(V(S)) = f'(V(S)) d_s V(S) \tag{60}$$

if $f$ is a smooth function on $R^d$.

$$d_t(V^1(S) d_s V(S)) = d_t V^1(S) d_s V(S) + V^1(S) d_t d_s V(S) \tag{61}$$

for two reals semi-martingales $V^1(S)$ and $V^2(S)$.

$$d_t(d_s V^1(S) d_s V^2(S)) = d_t d_s V^1(S) d_s V(S) + d_s V^1(S) d_t d_s V^2(S) \tag{62}$$

(For the various right-brackets for two-parameters semi-martingales involved in (62), we refer to [110].)

But, if we want that the semi-martingale lives on the manifold, we have to look at a more general class of semi-martingales than (33). Let us suppose that:

$$\begin{aligned} V(S) &= x + \int_{[0,S]} h(S') \delta B(S') + \int_{[0,S]} h^1(S') dS' \\ &+ \int_{[0,S]} h^2(S') \delta_s B(S') \delta_t B(S') + \int_{[0,S]} h^3(S') ds \delta_t B(S') \\ &+ \int_{[0,S]} h^4(S') dt \delta_s B(S') \end{aligned} \tag{63}$$

We consider Itô integrals in (63). Let us stress the difference between the stochastic integral $\int_{[0,S]} h(S') \delta B(S')$ and $\int_{[0,S]} h(S') \delta_s B(S') \delta_t B(S')$ which is quadratic in the leading Brownian sheet. (60), (61) and (62) remains formally valid for this class of semi-martingales. Let us suppose that $V(S)$ takes its values in $M$.

Since the Stratonovitch Calculus is the same than the traditional one, by Malliavin's transfer principle, $d_t V(S)$ can be seen formally as a process over $V(S)$ in the tangent bundle. So we cannot define intrinsically $d_s d_t V(S)$. We have to follow the requirement (33) in order to define $\nabla_s d_t V(S)$.

If $\omega$ is 1-form on the manifold, we get by Malliavin's transfer principle (see Section 7 about this principle):

$$d_s \omega(d_t V(S)) = \nabla \omega(d_s V(S), d_t V(S)) + \omega(\nabla_s d_t V(S)) \tag{64}$$

Let us recall in order to understand (64) that $\nabla \omega$ is defined by the formula for two vector fields $X$ and $Y$

$$\nabla \omega(X, Y) = X.\omega(Y) - \omega(\nabla_X Y) \tag{65}$$

where we take the derivative along the vector field $X$ of the function $\omega(Y)$.

Of course, (64) has to be seen at the level of Stratonovitch differential of two parameter semi-martingales.

Let us recall (see [110], p. 317): there exists a unique connection $\hat{\nabla}$ such that:

$$\nabla_s d_t V(S) = \hat{\nabla}_t d_s V(S) \tag{66}$$

$\hat{\nabla}$ preserves the Riemannian metric. r Let $\alpha(x, u, v)$ a bundle morphism from the principal bundle $O(M) \oplus O(M)$ over $M$ into the bundle $T(M) \otimes (R^m)^*$ over $M$. Let $B(S)$ be a $R^m$-valued Brownian sheet. The main theorem of Norris ([110], p. 318) is the following:

**Theorem 6** *The system of equations*

$$\nabla_s d_t x(S) = \alpha(d_s d_t B(S)) \tag{67}$$

$$\nabla_s u(s,t) = 0 \tag{68}$$

$$\hat{\nabla}_t v(s,t) = 0 \tag{69}$$

*has a unique solution if $x(0) = x$, $u(t,0) = u(0,s) = v(t,0) = v(0,s) = u$ for a semi-martingale $x(S)$ with values in $M$ and two semi-martingales $u(S)$ and $v(S)$ over $x(S)$ in $O(M)$.*

## 3. Ornstein-Uhlenbeck processes on loop spaces and Dirichlet forms

We are concerned in this part by a world-sheet equal to the cylinder $[0,1] \times S^1$, that is by diffusion processes on the loop space.

Let us consider the Brownian sheet $(t, s) \to B(t, s)$. It can be considered as a infinite dimensional Brownian motion $t \to (s \to B(t, s))$. Let $H_1$ be the Hilbert space of maps $h$ from $[0,1]$ into $R$ endowed with the energy norm $\int_0^1 |d/ds h(s)|^2 ds$. The Brownian sheet can be seen formally as the Brownian motion with values in the Hilbert space $H_1$. Let $h_i$ be an orthonormal basis of $H_1$. $H_1$ is densely continuously imbedded in $C([0,1], R)$ the space of continuous function from $[0,1]$ into $R$ endowed with the uniform norm. Let $F$ be a Fréchet smooth function on $C([0,1], R)$.

$$\Delta_G F = \sum F"(h_i, h_i) \tag{70}$$

exists and is finite. It is called the Gross Laplacian. Let us remark that the series in (70) converges because we consider an orthonormal basis of $H_1$ and because $F$ is Fréchet smooth on $C([0,1], R)$. The Brownian sheet is a Markov process on the path space governed by the Gross Laplacian.

There is another Laplacian on $C([0,1], R)$. It is the Ornstein-Uhlenbeck operator. Let us endow $C([0,1], R)$ with the Wiener measure $P$. Let $h$ in $H_1$. We have the following integration by parts formula:

$$E[<dF, h>] = E[F \operatorname{div} h] \tag{71}$$

where $\operatorname{div} h = \int_0^1 d/ds h(s) \delta B(s)$ if $F$ is a cylindrical functional

$$F(B(.)) = F^n(B(s_1), \ldots, B(s_n)) \tag{72}$$

and $B$ is a one dimensional Brownian motion and $F^n$ a smooth function from $R^n$ into $R$ with bounded derivatives of all orders. We define

$$<dF, h> = \sum \frac{\partial}{\partial x_i} F^n(B(s_1), \ldots, B(s_n)) h(s_i) \tag{73}$$

such that $dF$ appears as a random element of the dual of $H_1$. For a cylindrical functional

$$dF = \sum \frac{\partial}{\partial x_i} F^n(B(s_1), \ldots, B(s_n)) 1_{[0,s_i]}(s) \tag{74}$$

$dF$ is called the $H$-derivative, which is the basic tool of the Malliavin Calculus [112]. Since $H_1$ is a Hilbert space, we can assimilate $dF$ to a random element $\operatorname{grad} F$ of $H_1$. We put

$$\Delta_{O.H} = \operatorname{div} \operatorname{grad} \tag{75}$$

The definition is analog to the definition of the Laplace-Beltrami operator in finite dimension (see (3)), the Lebesgue measure which does not exist in infinite dimension being replaced by the Wiener measure.

These operations have an algebraic counterpart. Let $H_1^{sym,n}$ the $n^{th}$ symmetric tensor product associated to $H_1$, endowed with its natural Hilbert norm. Let $\Lambda$ the symmetric Fock space associated to $H_1$. To $\sigma = \sum h_n$ belonging to $\Lambda$, we associate the Wiener chaoses:

$$H(\sigma) = \sum \int_{[0,1]^n} h_n(s_1, \ldots, s_n) \delta B(s_1) \ldots \delta B(s_n) \tag{76}$$

(An element of the $n^{th}$ symmetric tensor product of $H_1$ can be assimilated namely as a symmetric application from $[0,1]^n$ into $R$, of finite $L^2$ norm. The map Wiener chaoses realize an isometry between the symmetric Fock space and $L^2(P)$. Under this identification,

$$< dF, H > = A_h F \tag{77}$$

where $A_h$ is the annihilation operator on the symmetric Fock space associated to $h$. (71) says nothing else that the adjoint of an annihilation operator on the symmetric Fock space is a creation operator and the Ornstein-Uhlenbeck operator is nothing else than the number operator which counts the length of the considered tensor product in the symmetric Fock space (see [104] for an extensive study).

The Ornstein-Uhlenbeck operator (75) is a symmetric positive self-adjoint densely defined operator on $L^2(P)$. It generates therefore a semi-group called the Ornstein-Uhlenbeck semi-group on the Wiener space. Moreover:

$$E[\Delta_{O,H} F.G] = E[< \operatorname{grad} F, \operatorname{grad} G >] \tag{78}$$

The right-hand side is called a Dirichlet form. Dirichlet forms are more tractable than Laplacians. There exists an abstract theory of Dirichlet forms and of process related. We refer to the course of Albeverio at Saint-Flour about that [3].

Let $E$ be a topological space, $\mu$ be a $\sigma$-finite positive measure on $E$. Let $C$ be a positive, symmetric, densely defined, closed bilinear form on $L^2(\mu)$. Closed means that if $F_n \to F$ in $L^2(\mu)$ and is a Cauchy sequence for $C$, $F$ belongs to the domain of $C$ and $C(F^n - F) \to 0$.

**Definition 7** $C$ is called a Dirichlet form if

$$C(\Phi(F)) \leq C(F) \tag{79}$$

for $\Phi(t) = 0$ if $t < 0$, $\Phi(t) = t$ if $t \in [0,1]$ and $\Phi(t) = 1$ if $t > 1$.

**Remark:** In fact, we have to regularize $\Phi$ (see [3], p. 36).

In general, Dirichlet forms are defined over a dense set and after, via some integration by parts, we perform the closure. (In (78), we consider cylindrical function, and after we perform the completion in order to get a closed Dirichlet form: it is possible, because we have integration by parts (71).) This procedure is an infinite dimensional generalization of the way to get Sobolev spaces in finite dimension, the Lebesgue measure being replaced by the Wiener measure. Namely the historical example of a Dirichlet form is when we take $E = R^d$ endowed with the Lebesgue measure.

$$C(F) = \int_{R^d} \sum_i (\frac{\partial}{\partial x_i} F)^2 dx \tag{80}$$

The Dirichlet form is closable on $R^d$ because we have integration by parts. The operator which is associated by the analog of (78) in this situation is nothing else than the Laplacian on $R^d$.

**Definition 8** Let $C$ be a Dirichlet form on $L^2(\mu)$. Let $O$ be an open subset of $E$. We define the capacity of $O$ as followed:

$$\operatorname{Cap} O = \inf_{F \in D(C)} (C(F) + \|F\|_{L^2(\mu)}^2) \tag{81}$$

where $F \geq 1$ almost surely on $O$. For any subset $A$ of $E$, we define

$$\operatorname{Cap} A = \inf_{A \subseteq O} \operatorname{Cap} O \tag{82}$$

for $O$ open subset of $E$.

A property is said to be satisfy quasi-everywhere if the property is satisfied outside a set of capacity 0.

**Definition 9** A Dirichlet form is said quasi-regular if:

(i) There exists compacts set $F_k$ such that $\operatorname{Cap}(F_k^c) \to 0$.
(ii) There exists a dense subset for $C$ of continuous functions on $E$, which separates the points of $E$.

**Remark:** The definition is in fact more general (see [3], Definition 24).

The basic result is therefore the following: a quasi-regular Dirichlet form defines outside a set of capacity 0 a process $t \to X_t(x)$. The process is continuous if the Dirichlet form is local:

$$C(F,G) = 0 \tag{83}$$

if the intersection of the support of $F$ and $G$ is of measure 0. The process is defined up a lifetime.

We will produce an example of quasi-regular local Dirichlet form on the based loop space, due to Driver-Roeckner [33].

Let be the heat semi-group $P_t = \exp[-t\Delta]$ on the compact Riemannian manifold $M$. It is represented by the Brownian motion and has a heat-kernel:

$$P_t f(x) = \int_M p_t(x,y) f(y) dg(y) \tag{84}$$

where $(x,y) \to p_t(x,y)$ is smooth strictly positive.

$dP^x$ is the law of the Brownian bridge starting from $x$ and coming back at $x$. $t \to \gamma(t)$ is a semi-martingale for $P^x$ where $t \to \gamma(t)$ is the canonical process on $L_x(M)$, the based loop space of $M$, that is the space of continuous maps $\gamma$ from $S^1$ into $M$ such that $\gamma(1) = x$.

The law $P^x$ is characterized by the following property: let $F^n(\gamma(s_1), \ldots, \gamma(s_n))$ $s_1 < s_2 < \cdots < s_n < 1$ be a cylindrical functional.

$$E[F^n] = p_1(x,x)^{-1} \int_{M^n} p_{s_1}(x, x_1) p_{s_2-s_1}(x_1, x_2) \ldots p_{1-s_n}(x_n, x) F(x_1, \ldots, x_n) dg(x_1) \ldots dg(x_n) \tag{85}$$

Let $\nabla$ be the Levi-Civita connection on $M$. We can define the parallel transport $\tau_t$ along $\gamma_t$ for the Levi-Civita connection. In local coordinates, it is the solution of the linear Stratonovitch differential equation:

$$d\tau_t = -\Gamma_{d\gamma_t} \tau_t \tag{86}$$

where $\Gamma$ denotes the Christoffel symbols of the connection (see [35] for a complete theory of horizontal lifts of semi-martingales).

$\tau_t$ realizes an isometry from $T_{\gamma(0)}(M)$ into $T_{\gamma(t)}(M)$. We can describe the space where the transport parallel lives in another way: let $e_t$ be the evaluation map $\gamma \to \gamma(t)$. It is a map from $L_x(M)$ into $M$. We introduce the pull-back bundle $e_t^* T(M)$ such that $\tau_t$ appears as a section of $(e_0^* T(M))^* \otimes e_t^* T(M)$.

The tangent space of differential geometry of a continuous loop $\gamma$ is realized by the set of continuous sections over $S^1$ of the bundle on the circle $\gamma^* T(M)$. But this natural tangent bundle of the loop space does not allow to do analysis.

In order to do analysis, we have to use the tangent bundle of Jones-Léandre [65] given in a preliminary form by Bismut [17]. We consider a section of $\gamma^* T(M)$ of the shape $\tau_t h(t) = X(t)$ where $h(.)$ belongs to $H_1$, endowed with the energy Hilbert structure $\int_{S^1} |d/ds h(s)|^2 ds = \|h\|^2 < \infty$, with boundary conditions $h(0) = h(1) = 0$. The tangent space of a loop is therefore an Hilbert space. If $h(.)$ is deterministic, we get Bismut's type integration by parts formula [31, 17, 56]

$$E[< dF, X >] = E[F \operatorname{div} X] \tag{87}$$

where

$$\operatorname{div} X = \int_0^1 < \nabla_s X(s), \delta\gamma(s) > + 1/2 \int_0^1 < S_{X(s)}, \delta\gamma(s) > \tag{88}$$

for a cylindrical functional given by a natural generalization to the manifold case by (72). $S$ in (88) is the Ricci tensor (see (117)) associated to the Levi-Civita connection and $\delta\gamma(s)$ the curved Itô integral on the manifold. It is therefore defined by $\int_0^1 <\nabla_s X(s), \delta\gamma(s)> = \int_0^1 <d/ds h(s), \delta B(s)>$ where $dB(s) = \tau_s^{-1} d\gamma(s)$.

$dF$ appears as a 1-form. Since the tangent space is a Hilbert space, $dF$ can be assimilated by duality as a measurable section $\operatorname{grad} F$ of the stochastic tangent bundle of $L_x(M)$. This point of view is intrinsic and does not require the study of a differential equation leading to the construction of the Brownian motion on $M$. The main theorem of Driver-Roeckner is the following:

$$C(F) = E[<\operatorname{grad} F, \operatorname{grad} F>] \tag{89}$$

defines a quasi-regular local Dirichlet form on the probability space $(L_x(M), P^x)$.

This means that we have to complete the Dirichlet form elementary defined for cylindrical functionals.

As an application of the abstract theory, we can define outside a set of capacity 0 a process $t \to X_t(\gamma)(s)$ on $L_x(M)$. Moreover, its lifetime is infinite. It is the Ornstein-Uhlenbeck process on the based loop space.

There are several extensions of the work of Driver-Roeckner.

For instance, Albeverio-Léandre-Roeckner [6] consider the free loop space $L(M)$, that is the set of continuous maps $\gamma$ from the circle into $M$. Albeverio-Léandre-Roeckner consider the Bismut-Hoegh-Krohn measure on $L(M)$:

$$d\mu = \frac{p_1(x,x) dg(x) \otimes dP^x}{\int_M p_1(x,x) dg(x)} \tag{90}$$

This measure is characterized as follows, if we neglect the normalizing term $\int_M p_1(x,x) dg(x)$, for a cylindrical functional $F = f_1(\gamma(s_1))\ldots f_n(\gamma(s_n))$ by:

$$E[F] = \operatorname{Tr}[\exp[-s_1 \Delta] f_1 \exp[-(s_2-s_1)\Delta] f_2 \ldots f_n \exp[-(1-s_n)\Delta]] \tag{91}$$

Albeverio-Léandre-Roeckner consider vector fields of the same type than before, but with the boundary conditions $\tau_1 h_1 = h_0$ (and not $h(0) = h(1) = 0$ because we consider the free loop space) with the Hilbert structure

$$\int_{S^1} \|h(s)\|^2 ds + \int_{S^1} \|d/ds h(s)\|^2 ds \tag{92}$$

By using the generalization of Bismut's type integration by parts formulas done by Léandre in [78] and [79], Albeverio-Léandre-Roeckner deduce a quasiregular local Dirichlet forms on the free loop space which is a natural extension of the Dirichlet form of (89). Moreover, over the free loop space, there is a natural circle action $\psi_t$: $\psi_t(\gamma)(s) = \gamma(s+t)$. Due to the cyclicity of the trace in (91), the Bismut-Hoegh-Krohn measure is invariant under the circle action (see [36] for a kind of reciproque) The Dirichlet form of Albeverio-Léandre-Roeckner is invariant under rotation. [6] deduce the existence outside a set of capacity 0 of a process, with lifetime infinite, invariant by rotation. It is the Ornstein-Uhlenbeck process on the free loop space.

Léandre [80] applies Dirichlet forms theory in order to define the Ornstein-Uhlenbeck process on the universal cover of the based loop space, when $\Pi^2(M)$ is different from 0.

Léandre-Roan [98] consider a developable orbifold $M/G$ where $G$ is a finite group. They use the description of Dixon-Harvey-Vafa-Witten [28] of the free loop space of the developable orbifold, which allows to avoid the difficulty to describe the singularities of the orbifold, in order to construct the Ornstein-Uhlenbeck process on the free loop space of the orbifold (An orbifold has singularities!).

## 4. Airault-Malliavin equation and infinite dimensional Brownian motion

We are concerned in this case by a world sheet which is $[0,1] \times N$ or $[0,1] \times S^1$ where $N$ is a compact manifold. We get a process $t \to \{S \to x(t,S)\}$ and the measure on the set of maps from $N$ into the target manifold $M$ $S \to x(1,S)$ is called following the terminology of Airault-Malliavin the heat-kernel measure. Namely, Airault-Malliavin [2] produce the theory of the Brownian motion on a loop group, which can be easily extended to the manifold case. Airault-Malliavin mechanism allows to produce random fields with parameter space any compact manifold into any compact manifold, with arbitrary regularity as it was extended by Léandre [82].

Let us recall that the theory of processes on infinite dimensional manifolds has a long history: see works of Kuo [74], Belopolskaya-Daletskii [13] and Daletskii [25]. Arnaudon-Paycha [10] have done too a theory of random processes on Hilbert manifolds.

We follow here the presentation of Léandre [82], which generalizes Airault-Malliavin equation to the case of any compact manifold as parameter space and to any compact manifold $M$ as target space.

Let $H_1$ be an Hilbert space continuously imbedded in the space of continuous functions from $N$ into $R$, endowed with the uniform topology. Let $S$ be the generic element of $N$ and $h$ be a generic element of $H_1$. We suppose:

**Hypothesis H:** There exists a map $(S, S') \to e_S(S')$ such that:

(i) $h(S) = <h, e_S>_{H_1}$
(ii) $(S, S') \to e_S(S')$ is Hoelder with Hoelder exponent $\alpha$.

We consider the Brownian motion $t \to B_t(.)$ with values in $H_1$.

$t \to B_t(S)$ is a R-valued finite dimensional Brownian motion and the right-bracket between $B_t(S)$ and $B_t(S')$ satisfies to:

$$<B_t(S), B_t(S')> = t e_S(S') \tag{93}$$

By Kolmogorov lemma [103], we deduce that $(t, S) \to B_t(S)$ has a version which is Hoelder.

As a matter of fact, we can consider an orthonormal basis $h_i$ of $H_1$ such that

$$B_t(S) = \sum B_i(t)h_i(S) \tag{94}$$

where $B_i(.)$ are some independent $R$-valued Brownian motions. We have

$$e_S(S') = \sum h_i(S)h_i(S') \tag{95}$$

But the series does not converge in $H_1$, but in a convenient space of Hoelder functions by **Hypothesis H**. We will describe the situation more precisely later.

We consider Airault-Malliavin equation

$$d_t x_t(x)(S) = \Pi(x_t(x)(S))d_t B_t(S); \quad x_0(x)(S) = x \tag{96}$$

where $B_t(.)$ is a collection of $m$ independent Brownian motion in $H_1$ still denoted $B_t(.)$.

This realizes a family of Brownian motion $t \to x_t(x)(S)$ in the manifold $M$. By using (93) and the Kolmogorov Lemma [103], we get ([82] Theorem 2.1):

**Theorem 10** $S \to x_1(S)$ *has almost surely a version which is* $\alpha/2 - \epsilon$ *Hoelder.*

**Scheme of the proof:** We have:

$$\begin{aligned} d_t x_t(x)(S) - d_t x_t(x)(S') &= (\Pi(x_t(x)(S) - \Pi(x_t(x)(S'))d_t B_t(S) \\ &\quad + \Pi(x_t(x)(S'))(d_t B_t(S) - d_t B_t(S')) \end{aligned} \tag{97}$$

¿From Burkholder-Davies-Gundy inequality and from (93) we deduce that

$$\begin{aligned} &E[|x_t(x)(S) - x_t(x)(S')|^p] \\ &\leq C \int_0^t E[|x_s(x)(S) - x_s(x)(S')|^p]ds + d(S,S')^{\alpha p/2} \end{aligned} \tag{98}$$

The result arises by Kolmogorov lemma [103] and Gronwall lemma.

**Remark:** [87] has considered the Sobolev space $H_k(N)$ of functions from $N$ into $R$

$$\int_N h(S)(\Delta_N + 1)^k h(S) dm(S) = \|h\|_k^2 \tag{99}$$

where $m$ is the Riemannian measure associated to a Riemannian structure on $N$ and $\Delta_N$ is the Laplace-Beltrami operator on $N$. The main property [51] is that $C^k$ norms of functions can estimated by their Sobolev norms in $H_{p(k)}$ where $p(k)$ tends to infinity when $k \to \infty$. Let $h_n$ be the sequence of eigenvectors of $\Delta_N + 1$ associated to the eigenvalues $\lambda_n$. When $n \to \infty$, $\lambda_n \sim Cn^\beta$ [51]. The Gaussian field associated to (99) can be written as

$$\sum \frac{1}{\lambda_n^k} B_n h_n = B(.) \tag{100}$$

where $B_n$ are independent centered normalized Gaussian variables. $E[\|h\|_{k'}^2] < \infty$ for some big $k'$ if $k$ is big enough. We apply Sobolev imbedding theorem [51].

Let $k_0$ be given. If $k$ is big enough, the Brownian motion with values in $H_k(N)$ is $C^{k_0+1}$ [100] and $S \to x_1(S)$ is almost surely of class $C^{k_0}$ due to the final dimensional Sobolev imbedding theorem, by an argument similar to the proof of Theorem 10.

Theorem 10 produces a generalization of Airault-Malliavin equation on a loop group: $N = S^1$, $M = G$ is a compact Lie group, $H_1$ is a convenient Sobolev space of maps from $S^1$ into the Lie algebra of $G$ and

$$d_t g_t(s) = g_t(s) d_t B_t(s) \qquad (101)$$

is the equation of the Brownian motion on the group $G$, when $s$ is fixed [2].

According to the terminology of Airault-Malliavin, the law of $S \to x_1(S)$ is called the heat-kernel measure on the set of maps from $N$ into $M$.

We would like to define the family of stochastic differential equations (96) as a process on the set of maps from $N$ into $M$. This requires a theory of random process on Banach manifolds, and, therefore, a theory of processes on Banach spaces. If we consider the set of maps from $N$ into $R$ endowed with some Hoelder norm, it is a very bad Banach space, because the Hoelder norm is very irregular. The good understanding of Stratonovitch differential equations on Banach spaces needs, following Brzezniak-Elworthy [19] the introduction of $M$-2 Banach manifolds. See [19] for the general theory.

We consider the space $H_1$ of maps $h$ from $S^1$ into $R^m$ endowed with the Hilbert structure:

$$\int_{S^1} |h(s)|^2 ds + \int_{S^1} |d/ds h(s)|^2 ds = \|h\|^2 \qquad (102)$$

Let $W^{\theta,p}$ be the sets of maps from $S^1$ into $R^m$ endowed with the Banach structure:

$$\|h\|_{\theta,p} = \left( \int_{S^1} |h(s)|^p ds + \int_{S^1 \times S^1} \frac{|h(s_1) - h(s_2)|^p}{|s_1 - s_2|^{1+\theta p}} ds_1 ds_2 \right)^{1/p} \qquad (103)$$

The interest of the norm $\|.\|_{\theta,p}$ is that it contains only integrals in its definition and not supremum as the Hoelder norm. In this point of view, the Sobolev-Slobodetski space $W^{\theta,p}$ is more interesting than an Hoelder Banach space. Moreover, if $1/p < \theta < 1$, the Brownian motion $B_t(.)$ with values in $H_1$ takes its values in $W^{\theta,p}$ [19]. The main ingredient of Brzezniak-Elworthy theory is:

**Theorem 11** *Let $F$ be the Nemytski map $h(.) \to (s \to \Pi(h(s)))$. Then $F$ is Fréchet smooth from $W^{\theta,p}$ into himself and has linear growth:*

$$\|F(h(.))\|_{\theta,p} \leq A + B \|h(.)\|_{\theta,p} \qquad (104)$$

The idea of Brzezniak-Elworthy is to consider the family indexed by $S^1$ of Stratonovitch differential equations on $M$

$$dx_t(s) = \Pi(x_t(s)) d_t B_t(s) \qquad (105)$$

as a unique Stratonovitch differential equation on $W^{\theta,p}$:

$$dX_t(.) = F(X_t(.))d_tB_t(.) \tag{106}$$

If there exists something like a nice martingale theory on $W^{\theta,p}$, we could prove, by reproducing the arguments of the theory of Stratonovitch differential equation on a finite dimensional manifold, by using Theorem 11, the following theorem:

**Theorem 12** *(106) defines a Markov process on $W^{\theta,p}$ if $p \in [2,\infty[, \theta \in ]0,1[$ and $1/p < \theta < 1$. Moreover, if the starting loops $\gamma$ belongs to $L(M)$, the free loop space of $M$, $t \to X_t(\gamma)$ belongs to $L(M)$.*

In order to prove this theorem, Brzezniak-Elworthy use the fact that $W^{\theta,p}$ is a $M$-2 Banach space if $p \in [2,\infty[$ and $\theta \in ]0,1[$. Let us recall what is a $M$-2 Banach space $E$.

A finite process indexed by the positive integers $M_k$ is called a martingale with values in $E$ respectively to the filtration $F_k$ if

$$E[M_{k'}|F_k] = M_k \tag{107}$$

if $k' > k$. We won't describe the natural integrability conditions which appear in the definition of a $E$-valued martingale (for instance $E[|M_k|^p] < \infty$). But, there is one which is satisfied for only special Banach spaces:

**Definition 13** *A Banach space is called M-2 if there exists a constant $C_2(E)$ such that for all martingales $M_k$:*

$$\sup_k E[|M_k|^2] \leq C_2(E) \sum_k E[|M_k - M_{k-1}|^2] \tag{108}$$

Let $\tilde{H}$ be an Hilbert space. We consider the formal Gaussian measure:

$$\frac{1}{\tilde{Z}} \exp(-\|\tilde{h}\|^2/2)dD(\tilde{h}) \tag{109}$$

with the formal Lebesgue measure $dD(\tilde{h})$. We suppose, as we have done several times, that there exists an inclusion $i$ from $\tilde{H}$ into $E$ such that $i(\tilde{H})$ is dense in $E$ such that the Gaussian measure lives on $E$. We say that we are in presence of an abstract Wiener space [75]. Let us consider the Brownian motion $\tilde{B}_t$ with values in $\tilde{H}$. It takes in fact its values in $E$. Let $\xi$ be a continuous linear map from $E$ into $E$. We deduce $\xi \circ i$ a continuous linear map from $\tilde{H}$ into $E$.

Let $\xi(t_k)$ be random continuous maps from $E$ into $E$, $F_{t_k}$ measurable for the $\sigma$-algebra spanned by the Brownian motion $\tilde{B}_t$ for $t \leq t_k$. We can define the elementary Itô integral with values in $E$:

$$I(x) = \sum_{k=0}^{n-1} \xi(t_k)(\tilde{B}_{t_{k+1}} - \tilde{B}_{t_k}) \tag{110}$$

¿From (108), since $t \to \tilde{B}_t$ is a martingale with values in $E$, we get:

$$E[|I(\xi)|^2] \leq C_2(E)E[\int_0^1 |\xi(t)|^2_{L(\tilde{H},E)}dt] \tag{111}$$

where $\xi(t) = \xi(t_k)$ if $t \in [t_k, t_{k+1}[$, a convenient subdivision of $[0,1]$.

Since we have the important property (111), Brzezniak-Elworthy can extend words by words the technics available for finite dimensional stochastic differential equation to the case of differential equations on Banach manifolds modelled on Banach spaces (in the Fréchet sense) of type $M-2$, if they consider **Fréchet-smooth** vector fields. (106) becomes a special case of this general theory.

We can characterize geometrically $M$-2 Banach spaces [115]. Pisier's characterization theorem allows to show that $W^{\theta,p}$ is a $M$-2 Banach space under the previous considerations.

Let us indicate the road of possible developments.

Maier-Neeb ( [101]) have defined a central extension by a finite dimensional commutative Lie group of a current group of maps from $N$ into $G$, a compact Lie group. It is possible to define of this central extension a structure of $M-2$ Banach Lie group modelled on a suitable Besov-Slobodetski space. We consider the Brownian Motion $B_t(.)$ in some Besov-Slobodetski space of maps from $N$ into the Lie algebra of $G$ and the Brownian motion $B^1(t)$(finite dimensional) in the Lie algebra of the finite dimensional commutative Lie group giving the central extension of the current group. We can apply Brzezniak-Elworthy theory in this situation in order to study the stochastic differential equation on the central extension of the current group:

$$d\hat{g}_t = \hat{g}_t d\hat{B}_t \tag{112}$$

where $\hat{B}_t = B_t(.) \oplus B_t^1$.

Léandre [84] consider a world sheet $\Sigma$ with a random Riemannian structure. In order to define the measure on the set of metrics $g$ on the world-sheet, Léandre uses Faddeev-Popov procedure, which allows to define the measure on the Teichmueller space, and the Shavgulidze measure on the gauge group [120]. This allows Léandre to produce a stochastic Polyakov measure, with a unitary action of the gauge group (the spaces of diffeomorphism of $\Sigma$ isotopic to the identity) on the $L^2$ of the theory (see [5] for a mathematical introduction to physicists bosonic string theory).

In fact, physicists are not only interested by Bosons, but too by Fermions. Let $\Lambda(R^n)$ be the exterior algebra of $R^n$. Let $\theta_I = \theta_{i_1} \wedge \cdots \wedge \theta_{i_I}$ if $I = (i_1, \ldots, i_I)$ where $\theta_i$ is the canonical basis of $R^n$. Let $F(\theta) = \sum a_I \theta_I$. The Berezin integral is defined by:

$$\int_B F(\theta) = a_{1,2,\ldots,n-1,n} \tag{113}$$

Physicists are interested by Fermionic Gaussian processes. Alvarez-Gaume [7], has deduced a super-symmetric proof of the Index theorem. A. Rogers [117] has defined a super-Brownian motion on a finite dimensional manifold which allows to give a rigorous version of the formal proof of Alvarez-Gaume of the Index theorem on a finite dimensional manifold. Witten [127] is motivated by the Index theorem on the free loop space of a manifold. He is motivated in fact by super-symmetric two-dimensional field theories. Motivated by [127], Léandre [86] has done an infinite dimensional extension of the work of A. Rogers and has defined the super-Brownian motion on a loop group.

## 5. Logarithmic Sobolev inequality and heat kernel measure

Let us consider a finite dimensional compact Riemannian manifold $M$. Let $dg$ be the Riemannian probability measure on $M$ and $\Delta_M$ be the Laplace-Beltrami operator on $M$. It has a discrete spectrum and can be diagonalized [51]. This implies that there is a spectral gap [26]. This means that:

$$\int_M f^2 dg - (\int_M f dg)^2 \leq C \int_M f \Delta_M f dg \tag{114}$$

for any function $f$ on $M$ for a suitable constant $C$. Spectral gap inequality is implied by Logarithmic Sobolev inequality [26]: If $\int_M f^2 dg = 1$

$$\int_M f^2 \log f^2 dg \leq C \int_M |\operatorname{grad} f|^2 dg \tag{115}$$

We had seen that Dirichlet forms can be defined in infinite dimension. In particular, it is known since a long time that the Dirichlet form (3.6) satisfies a Logarithmic Sobolev inequality. This implies that the Ornstein-Uhlenbeck operator on the Wiener space has a spectral gap (It is in this case obvious, because we can diagonalize the Ornstein-Uhlenbeck operator, which is too the Number operator on the Symmetric Fock space).

Let $dP_x$ be the law of the Brownian motion starting from $x$ in the compact manifold. $s \to \gamma_s$ denotes the canonical process on the based path space of paths starting from $x$ and $s \to \tau_s$ the stochastic parallel transport for the Levi-Civita connection called $\nabla$.

The curvature tensor is given by

$$R(X,Y)Z = \nabla_X \nabla_Y Z - \nabla_Y \nabla_X Z - \nabla_{[X,Y]} Z = R(X,Y)Z \tag{116}$$

for vector fields $X, Y, Z$ on $M$. It is tensorial in all its arguments. The Ricci curvature is given by:

$$S(Y) = \operatorname{Tr} R(.,Y)(.) = -\sum R(e_i, Y) e_i \tag{117}$$

where $e_i$ is an orthonormal basis of $T_x(M)$.

Bismut [17] considered the following differential equation

$$\tau_s^{-1} \nabla_s X_s^B(h) = d/ds\, h(s) ds - 1/2 \tau_s^{-1} S(X_s^B(h)) ds \tag{118}$$

where $h$ belongs to the Cameron-Martin space of maps from $[0,1]$ into $T_x(M)$ such that $h(0) = 0$ and such that $\int_0^1 |d/ds h(s)|^2 ds < \infty$. (We refer to (54) for the notation $\nabla_s X_s$.) If $X_s^B(h) = \tau_s H_s$, (118) means that $dH_s = dh_s - 1/2\tau_s^{-1} S(X_s^B(h)) ds$.

Let us explain the mysterious appearance of the Ricci tensor in (118) or in (88). We would like to compute $\nabla_{X_.} \tau_s^{-1} d\gamma_s$ where $X_t = \tau_t h_t$. For that we have to take the derivative of the parallel transport. $e_t^* T(M)$ appears as bundle on the path space, which inherits from a connection $\nabla^\infty$ from the Levi-Civita connection on the tangent bundle of $M$. If we use the considerations following (88), we see that:

$$\nabla_{X_.} \tau_t^{-1} d\gamma_t = \nabla_{X_.}^\infty \tau_t^{-1} d\gamma_t + \tau_t^{-1} \tau_t dh_t \tag{119}$$

By Bismut-Araf'evas formula [9, 17]

$$\nabla_{X_.}^\infty = \tau_t \int_0^t \tau_s^{-1} R(d\gamma_s, X_s) \tau_s \tag{120}$$

such that

$$\nabla_{X_.} \tau_t^{-1} d\gamma_t = dh_t - \int_0^t \tau_s^{-1} R(d\gamma_s, X_s) \tau_s dB_t \tag{121}$$

We have used the theory of Eells-Elworthy-Malliavin, which says that

$$d\gamma_s = \tau_s d_s B_s \tag{122}$$

where $dB_s$ is the Stratonovitch differential of a flat Brownian motion in $T_x(M)$.

The Ricci tensor appears when we convert the last Stratonovitch integral in an Itô integral. Tangent vectors of Jones-Léandre [65] are simple to write, but the Ricci tensor appears in the integration by part formula. For Bismut's tangent vector, vector fields are complicated to write but the integration by parts formula is simple to write. Namely, for a cylindrical functional, Bismut established the integration by parts formula [17]

$$\begin{aligned} E[<dF^B, X^B(h)>] &= E[<\mathrm{grad}^B F, X^B(h)>] \\ &= E[F^B \int_0^1 <\tau_s d/ds h(s), \delta\gamma(s)>] \\ &= E[F^B \int_0^1 <d/ds h(s), \delta B_S>] \end{aligned} \tag{123}$$

This allows Fang-Malliavin [38] to establish the following Clark-Ocone formula:

$$F(\gamma_.)] = E[F] = \int_0^1 <E[k_s|F_s], \delta B_s> \tag{124}$$

where $<dF^B, X_.^B(b)> = <k, h>_{H_1}$ and where $F_s$ is the filtration spanned by the curved Brownian motion. Let us explain in more details the difference between Bismut's way of differential Calculus on the path space and the way defined in Section 3 (see [79] for an extensive study of that).

We consider the equation

$$d\tau_t^B = -\Gamma_{d\gamma_t}\tau_t^B - 1/2S_{\tau_t^B}dt \tag{125}$$

starting from identity. A Bismut's vector field is written has $X_t^B(h) = \tau_t^B h_t$ with Hilbert metric

$$\|X^B(h)\|^2 = \int_0^1 |d/ds h_s|^2 ds \tag{126}$$

For a cylindrical functional $F = F^n(\gamma_{s_1}, \ldots, \gamma_{s_n})$

$$dF^B = \sum <\mathrm{grad}_{\gamma_{s_i}} F^n(\gamma_{s_1}, \ldots, \gamma_{s_n}), \tau_{s_i}^B 1_{[0,s_i]}(s).> \tag{127}$$

Clark-Ocone formula allows to Capitaine-Hsu-Ledoux [23] to show, by a simple use of the Itô formula, that, if $E[F^2] = 1$

$$E[F^2 \log F^2] \leq CE[\|\mathrm{grad}^B F\|^2] \tag{128}$$

On the other hand, Léandre [78], p. 521, has remarked that the tangent space of Bismut of the Brownian loop is isomorphic to Jones-Léandre tangent space [65] and that the norm of the isomorphism as well as its inverse is bounded. This allows Capitaine-Hsu-Ledoux [23] to show the following theorem:

**Theorem 14** $E[F^2 \log F^2] \leq CE[\|\mathrm{grad}\, F\|^2]$ *if* $E[F^2] = 1$ *where we follow the notations of (78).*

This constitutes a Logarithmic Sobolev inequality on the continuous based path space, whose first proof was done by Hsu [55]. Léandre [81] has done an analogous inequality on $C^1$ paths, by using Capitaine-Hsu-Ledoux method.

For material about logarithmic Sobolev inequality, we refer to [8].

We are motivated by an infinite dimensional generalization of this proof. We consider the based loop group $L_e(G)$ of continuous maps from $S^1$ into $G$ starting from the identity.

According [2], the tangent space of a loop $g_\cdot$ is constituted of loops of the type $X_s(h) = g_s h_s$ where $h_\cdot$ is a finite energy loop in the Lie algebra of $G$ starting from 0. We take as Hilbert norm

$$\|X_\cdot(h)\|^2 = \int_{S^1} |d/ds h_s|^2 ds \tag{129}$$

We can compute easily the bracket of two vector fields and show it is still a vector field. Therefore, by an abstract argument, there is a Levi-Civita connection on $L_e(G)$. We can compute its curvature. Freed [41] has succeeded to construct the Ricci tensor associated to the Levi-Civita connection on the based loop group, because we are in one dimension. These ingredients allow to Fang [37] to repeat the proof of Capitaine-Hsu-Ledoux of logarithmic Sobolev inequalities, but for the Brownian motion on a loop group instead of the Brownian motion on a finite dimensional compact manifold. Fang consider the Brownian motion $B_t(.)$ with values in $H_1$, the space of based loop in the Lie algebra of $G$ given by the Hilbert structure (129) and consider Airault-Malliavin equation:

$$d_t g_t(s) = g_t(s) d_t B_t(s); \quad g_0(s) = e \tag{130}$$

Fang consider the heat-kernel measure on $L_e(G)$ given by the law of $s \to g_1(s)$. Fang states for the heat-kernel measure the following integration by parts formula valid for cylindrical functionals:

$$E[< dF, X(h) >] = E[F \operatorname{div} X(h)] \tag{131}$$

This allows to define a Dirichlet form $E[< \operatorname{grad} F^1, \operatorname{grad} F^2 >]$. Fang [37], by using an infinite dimensional generalization of the method of Capitaine-Hsu-Ledoux [23] established the following theorem:

**Theorem 15** *For the heat-kernel measure on the based loop-group, we have if $E[F^2] = 1$*

$$E[F^2 \log F^2] \leq CE[\|\operatorname{grad} F\|^2] \tag{132}$$

We refer to [32, 62] and [63] for various statements about Logarithmic-Sobolev inequalities for heat-kernel measure on loop groups.

## 6. Wentzel-Freidlin Estimates

Let us look at the situation of Theorem 10:

Let $\tilde{H}$ the Hilbert space of maps from $[0, 1]$ into $H_2 = H_1 \oplus \cdots \oplus H_1$ endowed with the Hilbert structure:

$$\int_0^1 \|d/dt h_t\|_{H_2}^2 = \|h\|^2 \tag{133}$$

if $h_t \in H_2$.

The Brownian motion $B_t(.)$ with values in $H_2$ has formally as law the Gaussian measure on $\tilde{H}$

$$\frac{1}{\tilde{Z}} \exp(-\|h\|^2/2) dD(h) \tag{134}$$

where $dD(h)$ is the formal Lebesgue measure on $\tilde{H}$. This statement is clarified by the following large deviations estimates. Let $\epsilon B_t(.)$ the Brownian motion driven by a small parameter $\epsilon$. Following (40), (41), let $A$ be a Borelian subset of the set of maps from $[0, 1] \times N$ into $R^m$ endowed with the uniform topology. Int $A$ denotes its interior for the uniform topology, and Clos $A$ its closure for the uniform topology. We get, analogously to (40), (41), when $\epsilon \to 0$:

$$- \inf_{h \in \operatorname{Int} A} \|h\|^2 \leq \underline{\lim} 2\epsilon^2 P\{\epsilon B_.(.) \in A\} \tag{135}$$

$$\overline{\lim} 2\epsilon^2 P\{\epsilon B_.(.) \in A\} \leq - \inf_{h \in \operatorname{Clos} A} \|h\|^2 \tag{136}$$

This statement has an abstract counterpart: let $(H, E)$ be an abstract Wiener space. Let $e$ be the generic element of $E$ and $d\mu$ the Gaussian measure on $E$ associated to this abstract Wiener space. Let $\|.\|$ be the Hilbert norm on $H$ and $|.|$ the Banach norm on $E$. Let $A$ be a Borelian subset of $E$, Int $A$ its interior for the metric $|.|$ and Clos $A$ its closure for the metric $|.|$. We have the analog of (135) and (136). When $\epsilon \to 0$

$$- \inf_{h \in \text{Int } A} \|h\|^2 \leq \underline{\lim} 2\epsilon^2 \mu\{\epsilon e \in A\} \tag{137}$$

$$\overline{\lim} 2\epsilon^2 \mu\{\epsilon e \in A\} \leq - \inf_{h \in \text{Clos } A} \|h\|^2 \tag{138}$$

This abstract theorem can be for instance applied to the free field, where $E$ can be chosen as a negative order Sobolev space.

We consider Airault-Malliavin equation submitted to a small parameter $\epsilon$:

$$dx_t(\epsilon)(S) = \epsilon \Pi(x_t(\epsilon)(S))dB_t(S); \quad x_0(\epsilon)(S) = x \tag{139}$$

We cannot apply directly the abstract previous abstract results, because the solution of (139) is only almost surely defined. But we get the following lemma, analogous to Lemma 2.:

**Lemma 16** *For all $a, R, \rho > 0$, there exists $\alpha, \epsilon_0$ such that:*

$$P\{\|x.(\epsilon)(.) - x.(h)(.)\|_\infty > \rho; \|\epsilon B.(.) - h.(.)\|_\infty < \alpha\} \leq \exp[-R/\epsilon^2] \tag{140}$$

*for $\epsilon \in ]0, \epsilon_0]$ and $h$ such that $\|h\|^2 \leq a$ where*

$$dx_t(h)(S) = \Pi(x_t(h)(S))d_t h_t(S); \quad x_0(h) = x \tag{141}$$

**Scheme of the proof:** It is based upon an unwritten idea of Kusuoka [76] in order to prove large deviation estimates for flows. By a result of [11] (see [77] for a simple proof), we get: For all $a, R, \rho > 0$, there exists $\alpha, \epsilon_0$ such that:

$$P\{\|x.(\epsilon)(S) - x.(h)(S)\|_\infty > \rho; \|\epsilon B.(.) - h.(.)\|_\infty < \alpha\} \leq \exp[-R/\epsilon^2] \tag{142}$$

for $\epsilon \in ]0, \epsilon_0]$ and $h$ such that $\|h\|^2 \leq a$.

We define over $N$ a small triangulation where there are $\exp[C\epsilon^{-2}]$ sites $S_i$. We apply (142) to each site. We remark, by Kolmogorov lemma, that for some $\eta$:

$$P\{\sup_{t,S,S'} \frac{|x_t(\epsilon)(S) - x_t(\epsilon)(S')|}{d(S,S')^\eta} > R\exp[C'\epsilon^{-2}]\} \tag{143}$$

is smaller than $\exp[-K\epsilon^{-2}]$ for any $K$.

**Definition 17** *If $A$ is a Borelian subset of the continuous maps from $[0,1] \times N$ into $M$ endowed with the uniform topology, we put:*

$$\Lambda(A) = \inf_{x.(h)(.) \in A} \|h\|^2 \tag{144}$$

We get (see [82]):

**Theorem 18** *(Wentzel-Freidlin estimates).* *When $\epsilon \to 0$*

$$-\Lambda(\text{Int}\, A) \leq \underline{\lim} 2\epsilon^2 \log P\{x_.(\epsilon)(.) \in A\} \tag{145}$$

$$\overline{\lim} 2\epsilon^2 \log P\{x_.(\epsilon)(.) \in A\} \leq -\Lambda(\text{Clos}\, A) \tag{146}$$

For another proof, in the case of the Brownian motion on a loop group, we refer to the work of Fang-Zhang [39].

Moreover, in the last developments of string theory, [27], people look at branes: it is a submanifold $M_1$ of the target manifold $M$. If $\Sigma$ is a Riemann surface with boundary $\partial\Sigma$, people look at maps from $\Sigma$ into $M$ such that $\partial\Sigma$ is mapped on $M_1$. Léandre [88] consider the Brownian motion (in the sense of Airault-Malliavin) on the set of cylinders attached to a brane $M_1$ and performs the large deviation theory.

## 7. Stochastic Wess-Zumino-Novikov-Witten model and stochastic field theories

In fact, Felder-Gawdzki-Kupiainen [40] have shown that the Hilbert space $H$ of Segal's axioms should be the Hilbert space on $L^2$-sections of a convenient line bundle on the loop space, endowed with a convenient formal measure, and that the amplitudes of the theory should be connected with a generalized parallel transport of elements of this line bundle along the surface, which realizes the (interacting!) dynamics of the loop space.

Line bundles over loop spaces can be seen by two ways:

- Either the loop space is simply connected (this means that $\Pi^1(M) = \Pi^2(M) = 0$) and the line bundle is determined by its curvature.
- Either the loop space is not simply connected. There are constructions of line bundle of the loop space associated to Deligne cohomology of the manifold [45], gerbes [18] and bundle gerbes [49].

The line bundle on the loop space has to satisfy the so-called fusion property: let be two based loop in $x$ $\gamma^1$ and $\gamma^2$. We can construct, since the loop are based on the same point, a big loop $\gamma$. There is a map $\pi^1 : \gamma \to \gamma^1$ and $\pi^2 : \gamma \to \gamma^2$. We require that:

$$\xi(\gamma) = \pi^{1*}\xi(\gamma^1) \otimes \pi^{2*}\xi(\gamma^2) \tag{147}$$

Moreover, there are basically two theories:

- One on a compact Lie group. It is called the Wess-Zumino-Novikov-Witten model.
- One on a general compact manifold.

This part is relevant to the so-called **Malliavin's transfer principle:** A formula which is true in the deterministic context and which has a meaning via Stratonovitch stochastic Calculus, is still valid, but almost surely.

The philosophy of this part is to replace the kinetic term of conformal field theory by the more tractable one arising from Airault-Malliavin equation. The geometry is on the other hand very similar to the geometry developed in the surveys of Gawedzki [46, 47, 48]. Moreover, we can consider 1+2 dimensional theory, and the heat kernel measure associated or 1+1 dimensional theory.

We consider the torus $T^2 = S^1 \times S^1$. We consider the Hilbert space $H_1$ of maps $h(.)$ from $T^2$ into the Lie algebra of a compact Lie group

$$\|h\|^2 = \int_{T^2} <h(S), (-\frac{\partial^2}{\partial s^2} + 1)(-\frac{\partial^2}{\partial t^2} + 1)h(S) > dS \qquad (148)$$

if $S = (s, t)$.

This Hilbert satisfies to hypothesis H of Section 4. The Green kernel $E_S(S')$ is in this situation the product of the 1-dimensional Green kernels associated to the operator on the circle $\frac{-\partial^2}{\partial s^2} + 1$ which are equal to

$$e_0(s) = \lambda \exp[-s] + \mu \exp[s] \qquad (149)$$

for $0 \leq s \leq 1$ such that $e(0) = e(1)$. $e_s(s')$ is got by translation (see [85], p. 5534). Therefore $E_S(S')$ satisfies Hypothesis H of Section 4.

We consider Airault-Malliavin equation:

$$d_u g_u(S) = g_u(S) d_u B_u(S); \quad g_0(S) = e \qquad (150)$$

where $u \to B_u(.)$ is a Brownian motion with values in $H_1$. $S \to g_1(S)$ defines the heat-kernel measure $d\mu$ on the Hoelder torus group of Hoelder maps from $T^2$ into $G$ [85]. But, as a matter, of fact, $d\mu$ defines a measure on the strong Hoelder torus group $T^2_{\epsilon,*}(G)$ of maps such that:

$$\lim_{S \to S'} \frac{d(g(S), g(S'))}{d(S, S')^\epsilon} = 0 \qquad (151)$$

It is an infinite dimensional manifold. By the general theory of Bonic-Frampton-Tromba [14, 15], there are Fréchet smooth partition of unity associated to a cover of the strong Hoelder torus group by open subsets $V_i$ which is locally finite.

We get a generalization of Theorem 11 in this context [85]:

**Theorem 19** *Let $f$ be a map from $T^2 \times G$ into $G$ (conveniently imbedded into $R^m$) with bounded derivatives of all orders. Let $F$ be the Nemytski map:*

$$g(.) \to (S \to f(S, g(S))) \qquad (152)$$

*The Nemytski map is Fréchet smooth on the strong Hoelder torus group.*

This allows Léandre [85] to define stochastic plots, generalizing to this case the notion of diffeology of Chen [24] and Souriau [121] (see [43] for related works).

Let us recall what a diffeology on a topological space $M$ is. It is constituted of a collection of maps $(\phi_O, O)$ from any open subset $O$ of any $R^n$ satisfying the following requirements:

(i) If $j : O_1 \to O_2$ is a smooth map from $O_1$ into $O_2$ and if $(\phi_{O_2}, O_2)$ is a plot, $(\phi_{O_2} \circ j, O_1)$ is still a plot called the composite plot.
(ii) The constant map is a plot.
(iii) If $U_1$ and $U_2$ are two open disjoint subsets of $R^n$ and if $(\phi_{O_1}, O_1)$ and $(\phi_{O_2}, O_2)$ are two plots, the union map $\phi_{O_1 \cup O_2}$ realizes a plot from $O_1 \cup O_2$ into $M$.

This allows Chen and Souriau to define a form. A form $\sigma$ is given by data of forms $\phi_O^* \sigma$ on $O$ associated to each plot $(\phi_O, O)$. The system of forms over $U$ $\phi_U^* \sigma$ has moreover to satisfy the following requirement: if $(\phi_{O_2} \circ j, O_1)$ is a composite plot, $(\phi_{O_2} \circ j)^* \sigma$ is equal to $j^* \phi_{O_2}^* \sigma$.

Let us recall if $\sigma_O$ is a r-form on $O$ and if $j : O' \to O$ is smooth application $j^* \sigma(X_1', \ldots, X_r') = \sigma(DjX_1', \ldots, DjX_r')$ for vector fields on $O'$.

The exterior derivative $d\sigma$ of $\sigma$ is given by the data $d\phi_O^* \sigma$.

The main example of Souriau is the following: let $M$ be a manifold endowed with the equivalence relation $\sim$. We can consider the quotient space $\tilde{M}$. Let $\pi$ be the projection from $M$ onto $\tilde{M}$. A map $\tilde{\phi}$ from an open subset $U$ of a finite-dimensional linear space is a plot with values in $\tilde{M}$ if, by definition, there is a smooth lift $\phi$ from $U$ into $M$ such that $\tilde{\phi} = \pi \circ \phi$.

**Definition 20** *A stochastic plot of dimension $n$ on the strong Hoelder torus group is given by a countable family $(O, \phi_i, \Omega_i)$ where $O$ is an open subset of $R^n$ such that:*

(i) *The $\Omega_i$ constitute a measurable partition of $T^2_{\epsilon,*}(G)$.*
(ii) *$\phi(u)(S) = f_i(u, S, g(S))$ where $f_i$ is a smooth function over $O \times T^2 \times R^m$ with bounded derivatives of all orders into $R^m$.*
(iii) *Over $\Omega_i$, $\phi_i(u)(.)$ belongs to the torus group.*

*When there is no ambiguity, we call a stochastic plot $\phi$.*

We identify two stochastics plots $(O, \phi_i^1, \Omega_i^1)$ and $(O, \phi_j^2, \Omega_j^2)$ if $\phi_i^1 = \phi_j^2$ almost surely over $\Omega_i^1 \cap \Omega_j^2$.

This allows, following the lines of Chen [24] and Souriau [121] to define a stochastic form [85] (almost surely defined!):

**Definition 21** *A stochastic r-form $\sigma_{st}$ with values in $R$ is given by the following data: To each plot $\phi$ with source $O$, we associate a random smooth r-form $\phi^* \sigma_{st}$ on $O$ which satisfy to the following requirements:*

(i) *Let $j$ be a smooth map from $O^1$ into $O^2$ and let $\phi^2$ a plot with source $O^2$. We can consider the composite plot $\phi^1 = \phi^2 \circ j$ with source $O^1$. Then, almost surely:*

$$\phi^{1*} \sigma_{st} = j^* \phi^{2*} \sigma_{st} \tag{153}$$

(ii) *If $(O, \phi_i^1, \Omega_i^1)$ and $(O, \phi_j^2, \omega_j^2)$ are two stochastic plots such that $\phi_i^1 = \phi_j^2 \circ \psi$ on a set of probability different to 0 for a given measurable transformation of some $\Omega_i^1$ into some $\Omega_j^2$, then almost surely as smooth forms on $O$:*

$$\phi_i^{1*} \sigma_{st} = \phi_j^{2*} \sigma_{st} \circ \psi \tag{154}$$

We can define on the strong Hoelder torus group the standard de Rham cohomology groups, for forms smooth in Fréchet sense, because it is a Banach manifold. It is not heavy to do that, because there are partition of unity on the strong Hoelder torus group.

We can define the stochastic exterior derivative for stochastic forms $\sigma_{st}$. It is defined by the collection of $d\phi^*\sigma_{st}$ where we consider the collection of stochastic plots $\phi$. (See [24] and [121] in the deterministic context.)

[85] proves this theorem:

**Theorem 22** *The stochastic de Rham cohomology groups of the strong Hoelder torus group are equal to the deterministic Fréchet de Rham cohomology groups of the strong Hoelder torus group.*

The proof is based upon the fact there are partition of unity on the strong Hoelder torus group.

Let us recall quickly the material used in order to prove this theorem, which was used for a finite dimensional manifold in order to show that the various classical cohomologies of a finite dimensional manifold are equal [124].

Let $M$ be a topological space. It will be later the strong Hoelder torus group.

**Definition 23** *A presheaf $P = \{\Gamma(U; P); \rho(U, V)\}$ of R-vector spaces is a collection of vector spaces $\Gamma(U; P)$ indexed by the open subsets $U$ of $M$ and a restriction map $\rho(V; U) : \Gamma(V; P) \to \Gamma(U; P)$ for $U \subseteq V$ such that for $U \subseteq V \subseteq W$ we have:*
$$\rho(U; W) = \rho(U; V) \circ \rho(V; W) \tag{155}$$

**Definition 24** *A sheaf $\tilde{S} = \{\Gamma(U; \tilde{S}); \rho(U; V)\}$ of R-vector spaces is a presheaf such that for any cover $U_i$ by open subsets of $M$, the following two properties are checked:*

(i) *If the restriction to $U \cap U_i$ of a section $f$ belonging to $\Gamma(U; \tilde{S})$ equal the restriction to $U \cap U_i$ of another section $g$ of $\Gamma(U; \tilde{S})$, then $f = g$.*
(ii) *Let us give a system of section $f_i$ of $\Gamma(U_i; \tilde{S})$ such that the restriction to $U_i \cap U_j$ of $f_i$ is equal to the restriction to $U_i \cap U_j$ of $f_j$. There exists a unique $f \in \Gamma(U; \tilde{S})$ such that its restriction to $U_i \cap U$ are equal to the restriction of $f_i$ to $U \cap U_i$.*

If we replace $M$ by the strong Hoelder torus group, we can produce two sheaves on it:

(i) The sheaf of stochastic form in Chen-Souriau sense.
(ii) The sheaf of Fréchet smooth forms on it. It is possible to do that because the strong Hoelder torus group is a Fréchet manifold modelled on the space of strong Hoelder map from $T^2$ into the Lie algebra of $G$.

**Definition 25** *A morphism of sheaf $d : \tilde{S}' \to \tilde{S}$ is a collection of linear mappings $d_u$ from $\Gamma(U; \tilde{S}')$ into $\Gamma(U; \tilde{S})$ which are compatible with the operations of restrictions.*

For instance, a Fréchet smooth form on the strong Hoelder torus group defines clearly a stochastic form. We deduce a morphism of the sheaf of Fréchet smooth forms into the space of stochastic forms which is the inclusion.

The stochastic exterior derivative is a morphism of the sheaf of stochastic forms as well as the traditional exterior derivative for the sheaf of Fréchet smooth forms.

**Definition 26** *A morphism of sheaves $d : \tilde{S} \to \tilde{S}' \to \tilde{S}''$ is exact if for every open subset $V$, there exists an open subset $U \subseteq V$ such that $\text{Im}(d_U) = \text{Ker}(d_U)$.*

This means that we have a kind of Poincaré lemma. The stochastic exterior derivative is exact as well as the traditional exterior derivative on the sheaf of Fréchet smooth forms on the strong Hoelder torus group.

A sheaf over $M$ is said to be fine if for each locally finite cover $U_i$ by open subsets, there exists for each $i$ an endomorphism $l_i$ of the sheaf $S$ such that:

(i) $\text{Supp}\, l_i \subseteq U_i$.
(ii) $\sum l_i = 1$

$l_i$ are called partition of unity.

Since we work on the strong Hoelder torus group, the sheaf of stochastic form and of deterministic forms are fine.

Theorem 22 arises then by abstract arguments on sheaf cohomology [124].

Let us define a 1-dimensional stochastic plot $l$ with source $[a, b]$. We can define:

$$\int_l \sigma_{st} = \int_a^b l^* \sigma_{st} \tag{156}$$

if $\sigma_{st}$ is a stochastic 1-form (see [24, 121] in the deterministic case and [85] in the stochastic case).

We can consider a sum $l$ or a different of oriented segments $l^k$, with oriented boundaries. We say we are in presence of a stochastic 1-dimensional cycle if the boundaries destroy. We put:

$$\int_l \sigma_{st} = \sum \int_{l^k} \sigma_{st} \tag{157}$$

**Definition 27** *We say that a stochastic 1-form is Z-valued if $d\sigma_{st} = 0$ and if for any 1-dimensional stochastic cycle $l$, the random variable $\int_l \sigma_{st}$ belongs to $Z$ almost-surely.*

This definition is analogous to the deterministic one, with only difference that the stochastic form are involved, as we will see, with stochastic integrals, and therefore we cannot pick-up a trajectory of $g(.)$. We will perform the following hypothesis:

**Hypothesis H:** the torus group is connected and $\Pi^1(G) = 0$ such that the loop group is simply connected.

Namely if $\Pi(G) = 0$, $\Pi^2(G) = 0$.

In such a case, there exists a line (a 1-dimensional stochastic plot!) joining $e(.)$ to $g(.)$.

**Definition 28** *Let $\sigma_{st}$ be a Z-valued 1-form. Let $k$ be an integer. The generalized Wess-Zumino-Novikov-Witten model of level $k$ is given by the measure on the strong Hoelder torus group:*

$$d\mu_k = \exp[2\sqrt{-1}\pi k \int_l \sigma_{st}] d\mu \qquad (158)$$

This is defined independently of the connecting plot $l$, because $\sigma_{st}$ is Z-valued. This last property implies namely if two 1-dimensional stochastic plots $l_1$ and $l_2$ join $e(.)$ to $g(.)$, then $\int_{l_1} \sigma_{st}$ differs of $\int_{l_2} \sigma_{st}$ by a random integer.

This allows to define consistently the Wess-Zumino term $\exp[2\sqrt{-1}\pi k \int_l \sigma_{st}]$.

We will produce a stochastic interpretation of the classical examples of conformal field theory, by using Malliavin's transfer principle. This requires to define Stratonovitch integrals, since the field $g(.)$ is only Hoelder. It is done in [85] and is the object of the two next theorems:

**Theorem 29** *(Integral of a 1-form): Let $\omega$ be a 1-form on $G$. The stochastic integral of Stratonovitch type*

$$\int_{S^1} <\omega(g(s,t)), d_s g(s,t)> \qquad (159)$$

*exists almost surely and is limit in $L^2$ of the traditional random integrals $\int_{S^1} <\omega(g^N(s,t)), d_s g^N(s,t)>$ where $g^N$ is a convenient polygonal approximation of the random field $g(.)$.*

Let us show the mail estimate in order to show the existence of this stochastic integral (where we cannot apply Itô' theory of Stochastic integral, because there is no martingale involved with the definition of this (almost-surely defined!) stochastic integral). We remark is $s < s_1 + \Delta s_1 < s_2 < s_2 + \Delta s_2$ then

$$< B_.(s_1 + \Delta_{s_1}, t) - B_.(s_1, t), B_.(s_2 + \Delta s_2, t) - B_.(s_2, t) > = O(\Delta s_1 \Delta s_2) \quad (160)$$

By using Stochastic Calculus, we can show that:

$$\begin{aligned} E[&< f_1(g_1(s_1,t), g_1(s_1,t) - g_1(\Delta_{s_1} + s_1, t) > \\ &< f_2(g_1(s_2,t)), g_1(s_2,t) - g_1(s_2 + \Delta s_2, t) >] = O(\Delta s_1 \Delta s_2) \end{aligned} \qquad (161)$$

We compute the expectation of the square of polygonal approximation of the stochastic integral associated to (159) and we arrive to the sum expression like (161) which converges when the length of the associated subdivision tend to infinity.

**Theorem 30** *( Integral of a two-form): Let $\omega$ be a 2-form on $G$. The Stratonovitch stochastic integral $\int_{T^2} <\omega(g(S)), d_s g(S), d_t g(S)>$ exists and is limit in $L^2$ of the traditional random integrals $\int_{T^2} g^{N*}\omega$ where $g^N$ is a convenient polygonal approximation of the random field $g(.)$.*

The proof is based upon similar arguments than (160) and (161), but the combinatoric is much more complicated.

We consider the canonical 3-form $\omega$ on the simple simply connected Lie group $G$ which at the level of the Lie algebra of $G$ satisfies to:

$$\omega(X, Y, Z) = \frac{1}{8\pi^2} <[X, Y], Z> \qquad (162)$$

We have the following theorem [85]:

**Theorem 31**

$$\sigma_{st} = \int_{T^2} <\omega(g(S)), d_s g(S), d_t g(S), .> \qquad (163)$$

*defines a closed $Z$- valued stochastic 1-form on the strong Hoelder torus group.*

In order to see that, let us consider a stochastic plot $(\phi_O, O)$ given globally, in order to simplify, by $u \in O \to \{S \to f(u, S, g(S))\}$ (We take in order to simplify the exposition a global stochastic plot). Let $X$ be a vector field on $O$.

$$\phi^* \sigma_{st} X = \int_{T^2} <\omega(\phi(u)(S)), d_s \phi(u)(S), d_t \phi(u)(S), \partial_X \phi(u)(S)> \qquad (164)$$

The integral on a stochastic cycle in the torus group gives an integral on a random (continuous) 3-dimensional cycle on $G$. These cycle can be approximated by random piecewise smooth cycle in $G$: the integral of $\omega$ on these random smooth cycle is a random integer which converges by Theorem 30 to the integral of $\sigma_{st}$ on the stochastic one dimensional cycle of the torus group of $\sigma_{st}$.

The relation between the geometry of the torus group and the geometry of the loop groups comes from the following observation: from the heat-kernel measure on the torus groups, we deduced some measures on loop groups $L^t(G)$ given by $s \to g(s, t)$. We can define as before a stochastic diffeology on each loop group $L^t(G)$ as well as stochastic forms, stochastic two dimensional cycles with values in the loop group...

**Theorem 32**

$$\sigma^t_{st} = \int_{S^1} \omega = \int_{S^1} <\omega(g(S)), d_s g(s, t), ., .> \qquad (165)$$

*defines a Z-valued stochastic 2-form on each $L^t(G)$.*

Let $(\phi_O, O)$ be the above stochastic plot and $X, Y$ two vector fields on $O$. We have

$$\phi^* \sigma^t_{st}(X, Y) = \int_{S^1} <\omega(\phi(u(S)), d_s \phi(u)(S), \partial_X \phi(u)(S), \partial_Y \phi(u)(S)> \qquad (166)$$

The integral of $\sigma^t_{st}$ on a stochastic cycle in the loop group gives as above a random stochastic integral on a 3-dimensional cycle on $G$, which is a random integer as before.

Let $U_i$ be balls for the uniform distance on $L^t(G)$ of small radius and of center $g_i(., t)$. We can assume that they constitute a cover of $L^t(G)$. We can suppose that $\delta$ is small enough such that we can join $g_i(.,t)$ to $g(.,t) \in U_i$ by a stochastic segment, and therefore we can join the unit loop to $g(.,t)$ by a stochastic segment $l_i(g(.,t))$. We have supposed that the loop group is simply connected. So we can find a stochastic surface in $L^t(G)$ $B_{i,j}(g(.,t))$ whose boundary is $l_i(g(.,t))$ and $l_j(g(.,t))$ run in the opposite sense if $g(.,t) \in U_i \cap U_j$. We can define, analogously to (157)

$$\rho_{i,j}(g(.,t)) = \exp[-\sqrt{-1}2\pi k \int_{B_{i,j}(g(.,t))} \sigma_{st}^t] \qquad (167)$$

We get, since $\sigma_{st}^t$ is $Z$-valued:

– On $U_i \cap U_j$, almost surely:

$$\rho_{i,j}\rho_{j,i} = 1 \qquad (168)$$

– On $U_i \cap U_j \cap U_k$, almost surely:

$$\rho_{i,j}\rho_{j,k}\rho_{k,i} = 1 \qquad (169)$$

Namely $B_{i,j}(g(.,t)) \cup B_{j,k}(g(.,t)) \cup B_{k,i}(g(.,t))$ has no boundary such that the integral of $\sigma_{st}^t$ on it is a random integer.

This allows us to give the following definition:

**Definition 33** *The Hilbert space $H^t$ of $L^2$ section of the stochastic line bundle $\xi^t$ with curvature $2\pi k\sigma_{st}^t$ on $L^t(G)$ is constructed as follows. A section $\alpha^t$ of $\xi^t$ is given by a collection of random variable $\alpha_j^t$ on $U_j$ such that almost surely over $U_i \cap U_j$*

$$\alpha_j^t = \alpha_i^t \rho_{i,j}(g(.,t)) \qquad (170)$$

*Moreover, since $|\alpha_i^t| = |\alpha_j^t| = |\alpha^t|$, the Hilbert structure on $H^t$ is given by*

$$E[|\alpha^t|^2] = \|\alpha^t\|^2 \qquad (171)$$

A map from the torus into $G$ realizes a loop from $S^1$ into $L^1(G)$ by $t \to (s \to g(s,t))$. We can define a formal parallel $\tau$ transport along this loop for $\xi^1$ to $\xi^1$ and show that, almost surely:

$$\operatorname{Tr} \tau = \exp[2\sqrt{-1}\pi k \int_l \sigma_{st}] \qquad (172)$$

In order to see that, let us recall this basic statement [18]: let $M$ be a simply connected manifold. Let $\sigma$ be a closed 2-form, which is $Z$-valued. This means that the integral on any cycle of $\Sigma$ is an integer. $\sigma$ determines a unique line bundle with connection on $M$ whose curvature is $2\pi\sigma$. Let $l$ be a loop in $M$. Since $M$ is simply connected, $l$ is the boundary of a surface $\Sigma$. Then the holonomy of this line bundle for the given connection is given by:

$$\tau = \exp[2\pi\sqrt{-1}\int_\Sigma \sigma] \qquad (173)$$

We remark, since $\sigma$ is $Z$-valued that this expression does not depend of the chosen surface $\Sigma$ whose boundary is $l$.

This explains, following Gawedzki, the relation between the torus group and the loop group.

Léandre [89] has done a version to this theory to the case of the $1 + n$ punctured sphere: it is still a $1 + 2$ dimensional field theory. On the sphere, there is one input loop and $n$ output loops. Léandre adds some collars to the sphere. He get a heat-kernel measure for fields parametrized by $\Sigma(1 + n)$ which are only Hoelder. Moreover, the law of each loops at the boundary are the same and they are independents, because we have added these colors. By doing as in [85], this allows, by using stochastic integrals, to realize $\Sigma(1, n)$ as a map from $H^{\otimes n}$ into $H$ by using the fusion property (147). Moreover, the random fields parametrized by $\Sigma(1, n)$ are realized by sewing elementary pants. Moreover, there is a natural map

$$\Sigma(1, n) \times \Sigma(1, r_1) \times \cdots \times \Sigma(1, r_n) \to \Sigma(1, \sum r_i) \tag{174}$$

by sewing the exits loops of $\Sigma(1, n)$ on the input loops of $\Sigma(1, r_i)$. Léandre [89] showed, via the Markov property on the field on the sewing loops, that this operation of sewing punctured spheres along there boundaries is compatible with the natural composition maps:

$$\text{Hom}(H^{\otimes n}, H) \times \text{Hom}(H^{\otimes r_1}, H) \times \cdots \times \text{Hom}(H^{\otimes r_n}, H)$$
$$\to \text{Hom}(H^{\otimes \sum r_i}, H) \tag{175}$$

This statements say that the collection of $\text{Hom}(H^{\otimes n}, H)$ is an operad. Archetypes of operads are the set of trees. Relation between operads and two dimensional field theories was pioneered by Kimura-Stasheff-Voronov [68] and Huang-Lepowsky [58].

In [89], the geometry of $\Sigma(1, n)$ is fixed. Let us consider the case of a two dimensional field theory. $\Sigma$ is a Riemannian surface with exit and input boundary loops endowed with the canonical metric on $S^1$ on each on the connected components of the boundary [93].

Let us consider $\overline{\Sigma}$ got from $\Sigma$ by sewing disk along the boundaries. $\overline{\Sigma}$ has a canonical metric, inherited from $\Sigma$. Let $\Delta_{\overline{\Sigma}}$ be the Laplace Beltrami operator on $\overline{\Sigma}$. Let $H_{\overline{\Sigma}}$ the Hilbert space of maps $h$ from $\overline{\Sigma}$ into $R$ such that:

$$\int_{\overline{\Sigma}} (\Delta_{\overline{\Sigma}}^k + 1)h(S)(\Delta_{\overline{\Sigma}}^k + 1)h(S) dg_{\overline{\Sigma}}(S) < \infty \tag{176}$$

where $dg_{\overline{\Sigma}}(S)$ denotes the Riemannian measure on $\overline{\Sigma}$.

Let $B_{\overline{\Sigma}, t}$ be the Brownian motion with values in $H_{\overline{\Sigma}}$. If $k$ is big enough independent from $r$, $(t, S) \to B_{\overline{\Sigma}, t}(S)$ is continuous in $t \in [0, 1]$ and $C^r$ in $S \in \overline{\Sigma}$ (see [93]).

Let $\overline{S}_1 = [0, 1] \times S_1$ where we sew disk along the boundary. $\overline{S}_1$ inherits a canonical Riemannian structure. Let $H_{\overline{S}_1}$ be the Hilbert space of maps from $\overline{S}_1$ into $R$ such that

$$\int_{\overline{S}_1} (\Delta_{\overline{S}_1}^k + I)h(S)(\Delta_{\overline{S}_1}^k + 1)h(S) dm_{\overline{S}_1}(S) < \infty \tag{177}$$

Let $B_{\overline{S}_1,t}$ be the Brownian motion with values in $H_{\overline{S}_1}$. Let $g_\Sigma(S)$ be a map from $\Sigma$ into $[0,1]$ equal to 1 on $\Sigma$ where we have removed the output collars $[0,1/2[\times \Sigma_2$ and where we have removed the input collars $]1/2;1]\times \Sigma_1$. We suppose that $g_\Sigma$ is equal to zero on a neighborhood of the boundaries of $\Sigma$.

Let $g^{\text{out}}$ be a smooth map from $[0,1/2]$ into $[0,1]$ equal to 0 only in 0 and equal to 1 in a neighborhood of $1/2$. Let $g^{\text{in}}$ be a smooth map from $[1/2,1]$ equal to 0 only in 1 and equal to 1 in a neighborhood of $1/2$.

We consider the Gaussian random field parametrized by $\Sigma \times [0,1]$:

$$B_{\Sigma,.}(.) = g_\Sigma(.)B_{\overline{\Sigma},.}(.) + \sum_{\text{in}} g^{\text{in}} B^{\text{in}}_{\overline{S}_1,.}(.) + \sum_{\text{out}} g^{\text{out}} B^{\text{out}}_{\overline{S}_1,.}(.) \qquad (178)$$

where we take independent Brownian motion on $H_{\overline{S}_1}$ which are independent of the Brownian motion $B_{\overline{\Sigma}}$. We have a body process and some boundary processes which are independent themselves and of the body process.

An object $\Sigma_{\text{tot},k} = (\Sigma_1 \cup \Sigma_2 \ldots \cup \Sigma_k)$ is constructed inductively as follows: $\Sigma_1$ is a Riemann surface constructed as before. $\Sigma_{\text{tot},k+1}$ is constructed from $\Sigma_{\text{tot},k}$ where we sew some exit boundaries of $\Sigma_{\text{tot},k}$ along some input boundaries of $\Sigma_{k+1}$. Let us remark that in the present theory, we don't consider $\Sigma_{\text{tot},k}$ as a Riemannian manifold, but as the sequence $(\Sigma_1, \ldots, \Sigma_k)$ and the way we sew $\Sigma_{k+1}$ to $\Sigma_{\text{tot},k}$ inductively. Namely, if we consider the random fields parametrized by $\Sigma_{\text{tot},k} \times [0,1]$ considered as a global Riemannian manifold done by (176), it is different from the random field $B_{\Sigma_{\text{tot},k}}$ constructed as below. In particular, the sewing collars in $\Sigma_{\text{tot},k}$ are independent in the construction below, and are not independent in the construction (176).

We can construct inductively $B_{\Sigma_{\text{tot},k+1}}$ as follows: if $k = 1$, it is $B_\Sigma$. $B_{\Sigma_{k+1}}$ is constructed from Brownian motion independent of those which have constructed $B_{\Sigma_{\text{tot},k}}$, except for the Brownian motions in the input boundaries of $B_{\Sigma_{k+1}}$ which coincide with the Brownian motion in the output boundaries of $\Sigma_{\text{tot},k}$ which are sewed to the corresponding input boundaries of $\Sigma_{k+1}$. By this procedure, if $S \in \Sigma_{\text{tot}}$, we get a process $(t,S) \to B_{\Sigma_{\text{tot}},t}(S)$ which is continuous in $t$ and $C^r$ in $S \in \Sigma_{\text{tot}}$.

Let $G$ be the compact simply connected Lie group. We consider Airault-Malliavin equation [93]

$$d_t g_{\Sigma_{\text{tot}},t}(S) = g_{\Sigma_{\text{tot}},t}(S) \sum e_i d_t B^i_{\Sigma_{\text{tot}}}(S) \qquad (179)$$

starting from $e$. $B^i_{\Sigma_{\text{tot}}}$ are independent copies of $B_{\Sigma_{\text{tot}}}$ and $e_i$ an orthogonal basis of the Lie algebra of $G$.

**Theorem 34** *If $k$ is big enough, the random field parametrized by $\Sigma_{\text{tot}}$ $S \to g_{\Sigma_{\text{tot}},1}(S)$ is $C^r$. Moreover the restriction to this random field to the connected components of the boundary of $\Sigma_{\text{tot}}$ are independents and have the same law.*

Let us recall some geometrical background about the Wess-Zumino-Novikov-Witten model [42]. Let $\Sigma$ be an oriented surface with boundaries. Let $g$ be a $C^r$ map from $\Sigma$ into $G$ conveniently extended into a map $g_t(S)$ from $[0,1] \times \Sigma$ into

$G$ such that $g_0(z) = e$. We define the Wess-Zumino term:

$$W_\Sigma(g) = -1/6 \int_{[0,1]\times\Sigma} < g^{-1}dg \wedge [g^{-1}dg \wedge g^{-1}dg] > \qquad (180)$$

where $<,>$ is the canonical normalized Killing form on the Lie algebra of $G$. We suppose that the 3-form which is integrated in (180) represents an element of $H^3(G; Z)$ (see [42] for this hypothesis). $\exp[2\pi\sqrt{-1}W_\Sigma(g)]$ can be identified canonically to an element of $\xi_{\partial\Sigma,\partial g}$ where $\xi$ is an Hermitian line bundle over the set of $C^r$ maps from $\partial\Sigma$ into $G$. Let $\partial\Sigma_i$ be the oriented connected components of $\partial\Sigma$. We have a canonical inclusion map $\pi_i$ from $\partial\Sigma_i$ in $\partial\Sigma$. We deduce from it a map $\overline{\pi}_i$ from the set of maps from $\partial\Sigma$ in $G$ into the set of maps from $\partial\Sigma_i$ into $G$. Let $\xi_i$ be the hermitian bundle on the set of maps from $\Sigma_i$ into $G$ constructed in [42]. $\xi = \otimes\overline{\pi}_i^*\xi_i$ endowed with its natural metric inherited from each $\xi_i$. We denote it $\otimes_{\text{exit}}\xi\otimes_{\text{in}}\xi$. Moreover, we can realize this expression as a map from the tensor products of Hermitian line bundle $\xi$ over the exit loop groups defined by restricting the field over each exit boundary to the tensor product of Hermitian line bundles $\xi$ over the input loop groups defined by restricting the field over each connected component of the input boundary. Therefore $\exp[2\pi\sqrt{-1}W_\Sigma(g)]$ can be realized as an application from $\otimes_{\text{exit}}\xi$ into $\otimes_{\text{in}}\xi$ of modulus 1. This application is consistent with the operation of sewing surface.

Let $H'$ be the Hilbert space of section of $\xi$ over the $C^r$ loop group $L^r(G)$ endowed with the law arising from restricting the field to one boundary loops. Let $\Sigma_i$ be such a boundary loop. The laws of $g_{\Sigma_{\text{tot}},1}(.)$ restricted to each $\Sigma_i$ are the same. Let $\alpha_i(g_{\Sigma_{\text{tot}},1}(.)|_{\Sigma_i})$ a section of $\xi$ on the set of loops defined by $\Sigma_i$. We put

$$\|\alpha_i\|^2_{H_i'} = E[|\alpha_i(g_{\Sigma_{\text{tot}},1}(.)|_{\Sigma_i})|^2] \qquad (181)$$

Let $L^2([0,1]\times\Sigma_i)$ be the Hilbert space of $L^2$ functionals with respect of $g_{V_{\text{tot}},.}(.)$ restricted to $[0,1]\times\Sigma_i$. We put $H_i = H_i' \otimes L^2([0,1]\times V_i)$. We get always the same Hilbert space $H$ independent of the choice of $\Sigma_{\text{tot}}$.

**Definition 35** $A(\Sigma_{\text{tot}}, g_{\Sigma_{\text{tot}}})$ *is the operator from* $\otimes_{\text{out}}H$ *into* $\otimes_{\text{in}}H$ *where we put the tensor product along respectively the connected components of the exit boundary of $\Sigma_{\text{tot}}$ and of the input boundaries of $\Sigma_{\text{tot}}$ defined as follows: let $\alpha_i$ be a section of $\xi$ at the $i^{th}$ connected component of the exit boundary:*

$$A(\Sigma_{\text{tot}}, g_{\Sigma_{\text{tot}}}) \otimes_{\text{out}} \alpha_i$$
$$= E[\exp[2\pi\sqrt{-1}W_{\Sigma_{\text{tot}}}(g_{\Sigma_{\text{tot}},1})] \otimes_{\text{out}} \alpha_i | B'([0,1]\times\Sigma_1)] \qquad (182)$$

where $B'([0,1]\times\Sigma_1)$ is the $\sigma$-algebra spanned by the random field $g_{\Sigma_{\text{tot}},.}(.)$ restricted to the input data $[0,1]\times\Sigma_1$.

Let $(\Sigma^1_{\text{tot}}, g^1_{\Sigma_{\text{tot}}})$ and $(\Sigma^2_{\text{tot}}, g^2_{\Sigma_{\text{tot}}})$ and $(\Sigma_{\text{tot}}, g_{\Sigma_{\text{tot}}})$ got by sewing $\Sigma^1_{\text{tot}}$ along some exit boundaries coinciding with some input boundaries of $\Sigma^2_{\text{tot}}$. We call the sewing boundary $\tilde{\Sigma}$ in $\Sigma_{\text{tot}}$. By Markov property of the random field [71, 73, 105, 66, 116], we deduce [93]:

**Theorem 36** *We have:*

$$A(\Sigma_{\text{tot}}, g_{\Sigma_{\text{tot}}}) = A(\Sigma^1_{\text{tot}}, g_{\Sigma^1_{\text{tot}}}) \circ A(\Sigma^2_{\text{tot}}, g_{\Sigma^2_{\text{tot}}}) \tag{183}$$

*where the composition goes for the Hilbert spaces which arises from the sewing boundaries.*

Let us do a brief history of Markov property for random fields.

Let $(\Omega, F, P)$ be a probability space and $x(S)$ be a Gaussian continuous centered random field with parameter space a finite dimensional manifold $M$ equipped with a Riemannian distance $d$.

If $O$ is an open subset of $M$, we define

$$B(O) = \sigma(x(S); S \in O) \tag{184}$$

and for a closed subset $D$, we define:

$$B(D) = \cap_{\epsilon > 0} B(D_\epsilon) \tag{185}$$

where $D_\epsilon = \{S \in M, \inf_{S' \in D} d(S, S') < \epsilon\}$

**Definition 37** *A random field has the Markov property with respect to an open set $O$ if for all $B(\overline{O})$-measurable functional $F$:*

$$E[F|B(O^c)] = E[F|B(\partial O)] \tag{186}$$

*A random field is G-Markov if it has the Markov property with respect to all open set $O$.*

Markov property with respect to $O$ is equivalent to the following statement: for any event $A_1$ $B(\overline{O})$-measurable and for any event $A_2$ $B(O^c)$-measurable:

$$P(A_1 \cap A_2 | B(\partial O)] = P(A_1|B(\partial O)]P[A_2|B(\partial O)] \tag{187}$$

Let us recall that the reproducing Hilbert space $H$ of the continuous Gaussian random field is given as follows: if $F$ is a linear random variable of the Gaussian random field, we put:

$$f_F(S) = E[Fx(S)] \tag{188}$$

and

$$<f_F, f_{F'}> = E[FF'] \tag{189}$$

If $e_S(S')$ is the covariance of the continuous Gaussian random field:

$$E[x(S)x(S')] = e_S(S') \tag{190}$$

such that

$$h(S) = <h, e_S(.)> \tag{191}$$

($h$ is the generic element of the reproducing Hilbert space $H$ of the random field).

Let us recall ([73] Theorem 5.1):

**Theorem 38** *A random continuous Gaussian field $x(.)$ with reproducing Hilbert space $H$ is G-Markov if the two following conditions are checked:*

(i) *For all $h_1$ and $h_2$ with support disjoint $<f_1, f_2> = 0$.*
(ii) *If $H \in H$ is such that $f = f_1 + f_2$ with disjoint supports, then $f_1$ and $f_2$ belong to $H$.*

In particular, we have consider in this paper Gaussian G-Markov fields.

We were considering previously $1+2$ dimensional field theory, and the associated heat kernel measure and the case where the loop space is simply connected. We will consider now $1+1$ random fields (that is diffusion processes on loop spaces), with a non simply connected loop space. We recall briefly the construction of the Brownian pants of Brzezniak-Léandre [21] and Léandre [92].

We consider the compact Riemannian manifold $M$ of dimension $d$ isometrically imbedded in $R^m$. We introduce the Hilbert space $H_1$ of the set of loops in $R^m$ such that:

$$\int_0^1 |\gamma(s)|^2 ds + \int_0^1 |d/ds\gamma(s)|^2 ds = \|\gamma\|^2 < \infty \tag{192}$$

Moreover, the couple if $s \neq s'$ $t \to (B_t(s), B_t(s'))$ realizes a non degenerated Brownian motion over $R^m \times R^m$, although $t \to B_t(s)$ and $t \to B_t(s')$ are not independent: the covariance matrix of $B_t(s)$ and $B_t(s')$ is not degenerated. This comes from the fact that the two linear maps from $H_1$ into $R^m$ $\gamma(.) \to \gamma(s)$ and $\gamma(.) \to \gamma(s')$ are independents. The family of Stratonovitch equations

$$d_t x_t(s) = \Pi(x_t(s)) d_t B_t(s); \quad x_0(s) = x \tag{193}$$

has a meaning. We recall [92] that $(s,t) \to x_t(s)$ has almost surely a version which is $1/2 - \epsilon$ Hoelder for all $\epsilon$. This comes from the fact that the Green kernel of this theory are computed by (149) and satisfy Hypothesis H of Section 4.

Let $s_1 < s_2$ be two times. We constrain the elliptic diffusion $t \to (x_t(s_1), x_t(s_2))$ to be equal to $y$ at time 1.

Let us recall that if we consider an elliptic diffusion $\tilde{y}_t(\tilde{x})$ over a compact manifold $\tilde{M}$, it has an heat kernel $q_t(\tilde{x}, \tilde{y})$

$$E[f(\tilde{y}_t(\tilde{x}))] = \int_{\tilde{M}} q_t(\tilde{x}, \tilde{y}) f(\tilde{y}) dg(\tilde{y}) \tag{194}$$

satisfying the estimate:

$$|\operatorname{grad} \log q_t(\tilde{x}, \tilde{y})| \leq C \frac{\tilde{d}(\tilde{x}, \tilde{y})}{t} \tag{195}$$

for $\tilde{y}$ close to $\tilde{x}$ for the associated Riemannian metric and the natural Riemannian distance $\tilde{d}$ associated to the elliptic diffusion (see [17, 107]). Let us recall that if the stochastic differential equation of the elliptic diffusion is given by

$$d\tilde{y}_t(\tilde{x}) = \sum \tilde{X}_i(\tilde{y}_t(\tilde{x})) d\tilde{w}_t^i + \tilde{X}_0(\tilde{y}_t(\tilde{x})) dt \tag{196}$$

over the compact manifold, the bridge between $\tilde{x}$ and $\tilde{y}$ (that is the diffusion constrained in time 1 to be $\tilde{y}$ satisfies to the following stochastic differential equation (in Stratonovitch sense):

$$d\tilde{y}_t(\tilde{x},\tilde{y}) = \sum \tilde{X}_i(\tilde{y}_t(\tilde{x},\tilde{y})(d\tilde{w}_t^i + <\tilde{X}_i(\tilde{y}_t(\tilde{x},\tilde{y}),$$
$$\operatorname{grad}\log q_{1-t}(\tilde{y}_t(\tilde{x},\tilde{y}),\tilde{y}) > dt) + \tilde{X}_0(\tilde{y}_t(\tilde{x},\tilde{y}))dt \quad (197)$$

(see [17, 107]). This means that we transform $d\tilde{w}_t^i$ into $d\tilde{w}_t^i + \alpha_t^i dt$ by using the equation (196). By the estimate (195), we have:

$$E[\int_0^1 |\alpha_t^i| dt] < \infty \quad (198)$$

We write

$$B_t(s_2) = \alpha(s_1, s_2)B_t(s_1) + \beta(s_1, s_2)B_t(s_1, s_2)$$
$$B_t(s) = \alpha(s_1, s_2, s)B_t(s_1) + \beta(s_1, s_2, s)B_t(s_1, s_2) + \gamma(s_1, s_2, s)B_t(s_1, s_2, s) \quad (199)$$

where the Brownian motion $B_t(s_1, s_2, s)$ is independent of the Brownian motions $B_t(s_1)$ and $B_t(s_1, s_2)$. Conditionating by $x_1(s_1) = x_1(s_2) = y$ is nothing else to do the following transformation in (199).

$$d\tilde{B}_t(s) = \alpha(s_1, s_2, s)(dB_t(s_1) + \alpha_t^1(s_1, s_2)dt)$$
$$+ \beta(s_1, s_2, s)(dB_t(s_1, s_2) + \alpha_t^2(s_1, s_2)dt)$$
$$+ \gamma(s_1, s_2, s)dB_t(s_1, s_2, s) \quad (200)$$

We have:

**Lemma 39** *We can suppose (see (198)) that $\int_0^1 |\alpha_t^i| dt < K$ (**Hypothesis K**). Under this condition have:*

$$E[|x_t(s) - x_t(s')|^p]|x_1(s_1) = x_1(s_2) = y] \leq C|s - s'|^{p/2} \quad (201)$$

By using the Kolmogorov lemma (see [103]), we deduce that there exists an Hoelder version of the random field $(t, s) \to x_t(s)$ where we have conditionated by $x_1(s_1) = x_1(s_2) = y$. The loop $s \to x_1(s)$ is splitted in two loops $s \to x_1^1(s)$ and $s \to x_1^2(s)$ starting from $y$ and satisfying the estimates (201) if (**Hypothesis K**) is satisfied. Following the idea of Brzezniak-Léandre [21] and Léandre [92], we introduce two others Brownian motions with values in the Hilbert space $H_1$, $B_t^1(.)$ and $B_t^2(.)$, independent of each others and independent of the first Brownian motion $B_t(.)$. We consider the equations after time 1

$$d_t x_{t+1}^1(s) = \Pi(x_{t+1}^1(s))dB_t^1(s) \quad (202)$$

starting from the little loop $x_1^1(.)$ and

$$d_t x_{t+1}^2(s) = \Pi(x_{t+1}^2(s))dB_t^2(s) \quad (203)$$

starting from the second little loop $x_1^2(.)$. We have:

**Lemma 40** *If* (**Hypothesis K**) *is satisfied, we have:*

$$E[|x_t^1(s) - x_t^1(s')|^p] \leq C|s - s'|^{p/2} \qquad (204)$$

*and we have*

$$E[|x_t^2(s) - x_t^2(s')|^p] \leq C|s - s'|^{p/2} \qquad (205)$$

**Definition 41** *the random pant is constituted for $t \leq 1$ by the random field $(t, s) \to x_t(s)$ with the constrain $x_1(s_1) = x_1(s_2) = y$ and for $t > 1$ by the couple of diffusion processes $t \to (x_t^1(.), x_t^2(.))$.*

There are one input boundary at $t = 0$ and two output boundaries to the pant $(x_2^1(.), x_2^2(.))$.

Let us consider the product of loop spaces $L(M) \times L(M)$. We endow it with the probability law of $(x_2^1(.), x_2^2(.))$. We will construct a line bundle over $L(M) \times L(M)$, by using the arguments of Gawedzki [45]. We won't suppose that the loop space is simply connected, because our construction is motivated by Deligne cohomology [18]. Namely, in the case where the loop space is simply connected, the following construction are not useful because in such a case a (complex) line bundle is defined by its curvature, as we have seen before. The following constructions are interesting only when the loop space is not simply connected.

Let $O_\alpha$ be a cover of $M$ by convex contractibles open subsets of $M$, such that $O_{\alpha_1,\alpha_2} = O_{\alpha_1} \cap O_{\alpha_2}$, $O_{\alpha_1,\alpha_2,\alpha_3} = O_{\alpha_1} \cap O_{\alpha_2} \cap O_{\alpha_3}$ and $O_{\alpha_1,\alpha_2,\alpha_3,\alpha_4} = O_{\alpha_1} \cap O_{\alpha_2} \cap O_{\alpha_3} \cap O_{\alpha_4}$.

Let $g_{\alpha_1,\alpha_2,\alpha_3}$ be a family of smooth functions $S^1$-valued which are multiplicatively antisymmetric in $\alpha_1, \alpha_2, \alpha_3$ and such that

$$g_{\alpha_1,\alpha_2,\alpha_3} g_{\alpha_0,\alpha_2,\alpha_3}^{-1} g_{\alpha_0,\alpha_1,\alpha_3} g_{\alpha_0,\alpha_1,\alpha_2}^{-1} = 1 \qquad (206)$$

over $O_{\alpha_0,\alpha_1,\alpha_2,\alpha_3}$.

Also, let $\eta_{\alpha_1,\alpha_2} = -\eta_{\alpha_2,\alpha_1}$ be a smooth real 1-form over $O_{\alpha_1,\alpha_2}$ such that:

$$\eta_{\alpha_1,\alpha_2} - \eta_{\alpha_0,\alpha_2} + \eta_{\alpha_0,\alpha_1} = 1/i g_{\alpha_0,\alpha_1,\alpha_2}^{-1} dg_{\alpha_0,\alpha_1,\alpha_2} \qquad (207)$$

on $O_{\alpha_0,\alpha_1,\alpha_2}$. Finally, we suppose that $\omega_\alpha$ is a real 2-form defined on $O_\alpha$ such that:

$$\omega_{\alpha_1} - \omega_{\alpha_0} = d\eta_{\alpha_0,\alpha_1} \qquad (208)$$

on $O_{\alpha_0,\alpha_1}$. These data define an element of the second Deligne hypercohomology group of the manifold (see [18], p. 250–251). If we look at the 3-form $d\omega_\alpha = \omega$, they patch together by (208) in order to give a closed 3-form $\omega$ on $M$.

Consider a system $(l, v)$ which constitutes a triangulation of the circle $S^1$ such that $b$ is an edge and $v \in \partial b$ is one of its vertex. To each edge, we associate an element $\alpha_b$ and to each vertex $v$ we associate an number $\alpha_v$ such that the following hold: we consider the set of loops $\gamma$ such that for each edge $b$ $\gamma(b) \subseteq O_{\alpha_b}$ and such for all vertices $\gamma(v) \in O_{\alpha_v}$. This defines an open subset

$$U_{A,\alpha} = \{\gamma : S^1 \to M | \gamma(b) \subseteq O_{\alpha_b}, \gamma(v) \in O_{\alpha_v} \quad \text{for each } (b, v) \in A\} \qquad (209)$$

If we consider the product of the loop space, we consider the product $U_{A,\alpha} \times U_{A',\alpha'}$ which constitutes a cover by open subsets of the product of the loop space.

We would like to define a system of transition maps of $(U_{A_1\alpha_1} \times U_{A'_1,\alpha'_1}) \cap (U_{A_2,\alpha_2} \times U_{A'_2,\alpha'_2})$. Let us define the refined triangulation of both triangulation $A_1$ and $A_2$ by $(\bar{b},\bar{v})$, $\bar{v} \in \partial \bar{b}$, the triangulation $A_1$ by $(b_i, v_i)$, $v_i \in \partial b_i$ and the second triangulation by $(b_2, v_2)$, $v_2 \in \partial b_2$. Let us put $\alpha^1_{\bar{b}} = \alpha^1_{b_1}$ and $\alpha^2_{\bar{b}} = \alpha^2_{b_2}$. If $\bar{v}$ is a vertex of the new triangulation, we put $\alpha^1_{\bar{v}} = \alpha^1_{v_1}$ if $\bar{v} = v_1$ and $\alpha^1_{\bar{v}} = \alpha^1_{b_1}$ if $\bar{v}$ is in the interior point of the interval $b_1$. We define $\alpha^2_{\bar{v}}$ analogously. The system of transition functionals of the stochastic line bundle over $L(M) \times L(M)$ is defined by

$$\rho = \rho_{A_1,\alpha_1,A_2,\alpha_2}(x^1_2)\rho_{A'_1,\alpha'_1,A'_2,\alpha'_2}(x^2_2) \tag{210}$$

where

$$\rho_{A_1,\alpha_1,A_2,\alpha_2}(x^1_2) = \exp[i \sum_b \int_{\bar{b}} \eta_{\alpha^1_{\bar{b}}\alpha^2_{\bar{b}}}(d_s x^1_2(s))] \prod_{\bar{v},\bar{b},\bar{v} \in \partial\bar{b}} \frac{g_{\alpha^1_{\bar{v}}\alpha^2_{\bar{v}}\alpha^2_{\bar{b}}}}{g_{\alpha^1_{\bar{v}}\alpha^1_{\bar{b}}\alpha^2_{\bar{b}}}}(x^1_2(v)) \tag{211}$$

and the analogous formula holds for $\rho_{A'_1,\alpha'_1,A'_2,\alpha'_2}(x^2_2)$. The transition functions are almost surely defined, due to the presence of stochastic integral in the definition of them. So we cannot define $\xi = \xi_1 \otimes \xi_2$, but we will follow the lines of Definition 32 in order to define the Hilbert space of $L^2$ sections of it.

**Definition 42** *A $L^2$ section of the line bundle $\xi_1 \otimes \xi_2$ over $L(M) \times L(M)$ is a system of functionals over $U_{A,\alpha} \times U_{A',\alpha'}$ $\alpha_{A,\alpha,A,\alpha'}$ submitted to the relations: almost surely, over $(U_{A_1,\alpha_1} \times U_{A'_1 \times \alpha'_1}) \cap (U_{A_2,\alpha_2} \times U_{A'_2,\alpha'_2})$, we get $\alpha_{A_1,\alpha_1,A'_1,\alpha'_1} = \rho\alpha_{A_2,\alpha_2,A'_2,\alpha'_2}$.*

We can define since $\rho$ defined by (211) is of modulus 1 the norm of a section $|\alpha|$. We suppose $E[|\alpha|^2] < \infty$ in order to define the space of $L^2$ sections of $\xi_1 \otimes \xi_2$.

In order this definition has some consistency, we recall that almost surely over $(U_{A_1,\alpha_1} \times U_{A'_1,\alpha'_1}) \cap (U_{A_2,\alpha_2} \times U_{A'_2,\alpha'_2})$, we get:

$$\begin{aligned}\rho_{A_1,\alpha_1,A_2,\alpha_2}(x^1_2)\rho_{A_2,\alpha_2,A_1,\alpha_1}(x^1_2) &= 1 \\ \rho_{A'_1,\alpha'_1,A'_2,\alpha'_2}(x^2_2)\rho_{A'_2,\alpha'_2,A'_1,\alpha'_1}(x^2_2) &= 1\end{aligned} \tag{212}$$

and that on $U_{A_1,\alpha_1} \cap U_{A_2,\alpha_2} \cap U_{A_3,\alpha_3}$, we get almost surely:

$$\rho_{A_1,\alpha_1,A_2,\alpha_2}(x^1_2)\rho_{A_2,\alpha_2,A_3,\alpha_3}(x^1_2)\rho_{A_3,\alpha_3,A_1,\alpha_1}(x^1_2) = 1 \tag{213}$$

This identity still works for the product of transition functions defined by $\rho_{A_1,\alpha_1,A_2,\alpha_2}(x^1_2) \rho_{A'_1,\alpha'_1,A'_2,\alpha'_2}(x^1_2)$.

This comes from the fact these relations are surely true for the deterministic loop space and that we can approach in the stochastic case the stochastic integrals which appear in the (almost surely defined!) transition function by their polygonal approximation.

In the previous definition, we have supposed that the section is almost surely defined over the product of random loops $(x^1_2, x^2_2)$ and is $(x^1_2, x^2_2)$ measurable.

We can suppose that $\alpha_{A,\alpha,A',\alpha'}$ depends from all the random pants, or if we choose $1 \leq t \leq 2$, $1 \leq t \leq 2$, it depends from all the paths between $t$ and 2. $\alpha_{A,\alpha,A',\alpha'}$ becomes an element of $L^2(\text{pant}) \otimes L^2(U_{A,\alpha}(x_2^1) \otimes U_{A',\alpha'}(x_2^1))$ and still satisfies to the consistency relations of Definition 41. We can define the $L^2$ norm of the section $\alpha$. This increases the degree of freedom and is done in order to define what is the parallel transport over the random path from $(x_2^1, x_2^2)$ into the path $(x_1^1, x_1^2)$. We will get a section of the bundle over $L_x(M) \times L_x(M)$, $\xi_1 \otimes \xi_2$ for the measure defined by $(x_1^1, x_1^2)$, but with an extra degree of freedom, that is the path between $(x_2^1, x_2^2)$ to $(x_1^1, x_1^2)$. In order to define the stochastic parallel transport from a random section over $\xi_1 \otimes \xi_2$ over $(x_2^1, x_2^2)$ to $(x_1^1, x_1^2)$ along the path $t \to (x_t^1, x_t^2)$, we will use the double integral of the previous part.

Let us divide the time interval $[1, 2]$ into the stochastic intervals $[\tau_i, \tau_{i+1}[$ where $(x_t^1, x_t^2)$ over $[\tau_i, \tau_{i+1}[$ the process $(x_t^1, x_t^2)$ lives over some open subset $U_{A,\alpha} \times U_{A',\alpha'}$. We have described the time interval into a finite numbers of random intervals. Moreover, the times $\tau_i$ are stopping times.

Let us suppose that the parallel transport from $(x_{\tau_i}^1, x_{\tau_i}^2)$ to $(x_1^1, x_1^2)$ is well defined. Let us call it $\tau_{1,\tau_i} = \tau_{1,\tau_i}^1 \otimes \tau_{1,\tau_i}^2$ (The product formula will be explained by the next considerations). If $t \in [\tau_i, \tau_{i+1}[$, we have, by using analogous theorem in the present situation than Theorem 29 and Theorem 30,

$$\begin{aligned} \tau_{1,t} &= \tau_{1,\tau_i} \{\exp[\sqrt{-1} \int_{\tau_i}^{t \wedge \tau_{i+1}} \sum_b \omega_{\alpha_b}(d_s x_t^1(s), d_t x_t^1(s)) \\ &\quad + \sqrt{-1} \sum_{v,b,v \in \partial b} \int_{\tau_i}^{t \wedge \tau_{i+1}} \eta_{\alpha_v, \alpha_b}(d_t x_t^1(v))]\} \\ &\quad \otimes \{\exp[\sqrt{-1} \int_{\tau_i}^{t \wedge \tau_{i+1}} \sum_{b'} \omega_{\alpha_{b'}}(d_s x_t^2(s), d_t x_t^2(s)) \\ &\quad + \sqrt{-1} \int_{\tau_i}^{t \wedge \tau_{i+1}} \sum_{v',b',v' \in \partial b'} \eta_{\alpha_v, \alpha_b}(d_t x_t^2(v))]\} \end{aligned} \qquad (214)$$

Let us explain this formula: for the smooth loop space, [45] gives the formula in term of double integral of the amplitude of this line bundle (that is a generalized parallel transport). We can deduce in the system of local charts given in (209) the Connection 1-form of this line bundle on the loop space. If $\Gamma$ is a connection 1-form of a line bundle, the parallel transport is given by (57), which can be integrated since we consider a line bundle. The extension in the stochastic case in our situation gives (214).

Let us remark that by induction the parallel transport is of modulus one (214). The rules given in the previous parts of approximation of Stratonovitch integrals allow to state this theorem:

**Theorem 43** *If $\alpha$ is a section of $L(M) \times L(M)$ for the measure of $(x_2^1, x_2^2)$ and measurable for $(x_2^1, x_2^2)$, $\tau_{1,2}\alpha$ is a section of $\xi_1 \otimes \xi_2$ for the measure of $(x_1^1, x_1^2)$ (but in an extended sense, because there are many paths joining $(x_2^1, x_2^2)$ to $(x_1^1, x_1^2)$). Moreover,*

$$E[|\tau_{1,2}\alpha|^2] = E[|\alpha|^2] \qquad (215)$$

We refer to [92] for this result.

Let us work in time 1. We consider the product of loop space $L_y(M) \times L_y(M)$ for the measure $(x_1^1, x_1^2)$ and the loop $L(M)$ induced by concatenation of the two loops for the measure induced by $x_1$. This induces a map $\pi$

$$L_y(M) \times L_y(M) \to L(M) \qquad (216)$$

which preserves the measure. Over $L_y(M) \times L_y(M)$, we have the stochastic line bundle $\xi_1 \otimes \xi_2$ and over $L(M)$ we have the stochastic line bundle $\xi$ defined by the previous considerations for the random loop $x_1$.

For $x_1$, we define a triangulation by choosing vertices $s_1$ and $s_2$. We have $\gamma(s_1) = \gamma(s_2) = y \in O_{\alpha_y}$. We choose another triangulation $(b^1, v^1)$ where we have chosen $s_1$ and $s_2$ among the vertices. For the first triangulation, we suppose $\gamma(b) \subseteq O_{\alpha_b}$ and $\gamma(v) \in O_{\alpha_v}$ where $O_{\alpha_{s_1}} = O_{\alpha_{s_2}} = \overline{O}$ is fixed, and for the second triangulation, we choose $\gamma(b') \subseteq O_{\alpha_{b'}}$ and $\gamma(v') \in O_{\alpha_{v'}}$ where $O_{\alpha_{s_1}} = O_{\alpha_{s_2}} = \overline{O}$ for the same open subset $\overline{O}$ than the first triangulation.

We deduce from the previous triangulation two triangulations of $L_y(M)$ and from the second triangulation two triangulations of $L_y(M)$. The transition map for the big loop space $L(M)$ is given by (211) where we replace $d_s x_2^1(s)$ by $d_s x_1(s)$ and $x_2^1(v)$ by $x_1(v)$ for the refined triangulation of the two big triangulations of the big circle. But it is almost surely equal to the product of the two transition functions where we consider the couple of loops $(x_1^1(s), x_2^1(s))$. This shows us that the fusion property (147) $\xi_1 \otimes \xi_2 = \pi^* \xi$ is satisfied. This means that a $L^2$ section of $\xi$ over the random loops $x_1$ for $L(M)$ corresponds naturally to a $L^2$ section of $\xi_1 \otimes \xi_2$ over $L_y(M) \times L_y(M)$ for the law $(x_1^1, x_1^2)$ and the $L^2$ norms are conserved. We assimilate $\tau_{1,2}\alpha$ to a section over $x_1$. Afterwards, we use the stochastic parallel transport from $x_1$ to $x_0$ $\tau_{0,1}$. We put $\tilde{\alpha} = \tau_{0,1}\tau_{1,2}\alpha$. Since $x_0$ is the constant loop $s \to x$, $\tilde{\alpha}$ is a random variable. We have [92]:

**Theorem 44** $E[|\tilde{f}|^2] = E[|f|^2]$.

This comes from the fact that $\tau_{0,1}$ is a random isometry from the stochastic fiber of the bundle over the random loop $x_1$ to the fiber over the constant loop.

Brzezniak-Léandre [21] have used the theory of Brzezniak-Elworthy [19] in order to realize stochastic pants over a suitable Besov-Slobodetsky space $W^{\theta,p}$ of loops $\gamma(.)$ in the manifold $M$. The pants starts from the loop $\gamma(.)$ and has two end loops $x_2^1(\gamma)$ and $x_2^2(\gamma)$. Let $E$ be the Banach space of continuous bounded functionals on $W^{\theta,p}$ and $E \otimes E$ be the Banach space of bounded continuous functionals on $W^{\theta,p} \times W^{\theta,p}$. Let

$$T(F) : \gamma(.) \to E[F(x_2^1(\gamma), x_2^1(\gamma))] \qquad (217)$$

if $F \in E \otimes E$. $T$ realizes a continuous linear application from $E \otimes E$ into $E$. This means that the Brownian pant is Fellerian.

Brzezniak-Léandre [20] have applied the theory of Brzezniak-Elworthy [19] to some Besov-Slobodetsky space of differentiable loops, where the line bundle of Gawedzki [45] is surely defined, and have defined the stochastic parallel transport on it for the Brownian motion on differentiable loops.

Gawedzki-Reis [49] have a simpler way, with more tractable combinatorial formulas, to construct a line bundle on differentiable paths, by using bundle gerbes theory. Léandre [90] constructs the Brownian motion along the Hoelder loop space, defines a stochastic line bundle in the manner of Gawedzki-Reis and study the stochastic parallel transport of it along the path of the Brownian motion on the Hoelder loop space.

With this branching type mechanismus (classical in theoretical physic [102], Léandre [94] has defined a kind of Branching process on the loop space.

In (201), we have conditionated by $x_t(s) = x_t(s')$ for two fied time $s = s'$. When there is no cut-locus, we can conditionated by this procedure by $x_1(s) = y$ for all $s$, which produced a kind of Brownian bridge in infinite dimension ([L$_{19}$]!.

[95] produces a kind of conditionating by an infinite dimensional constrain, where we cannot apply the classical tool of Airault-Malliavin-Sugita construction.

By using tools of Airault-Malliavin-Sugita construction, Léandre [96] has produced a stochastic regularization of the Poisson-Sigma model of Cattaneo-Felder which gives an infinite dimensional analog of Klauder's regularization of Hamiltonian path integral in quantum mechanic [67]. So there are two regularization in field theory:

(i) The first one is stochastic quantization of Parisi-Wu, which uses infinite dimensional Langevin equation.
(ii) The second one is stochastic quantization of Klauder, which uses infinite dimensional Brownian motion of Airault-Malliavin.

Let us remark that Brylinski [18] constructed a line bundle by using category theory and gerbes theory of Grothendieck on the smooth loop space. Léandre [91] gives a stochastic interpretation of [18] by studying a stochastic line bundle on the Brownian bridge of a manifold.